	\newcommand{\open}{\Bbb}
	 
	\newlength{\labparwidth}
	\newcommand{\se}{\subseteq}

	\newcommand{\gL}{\Lambda}

	\documentstyle[12pt,amsfonts,titlepage]{article}

	\title{A Variety with Solvable, but not Uniformly Solvable, Word
	Problem} 
	\author{Alan H. Mekler\\Department of Mathematics and
	Statistics\\Simon Fraser University\\Burnaby, B.C. V5A 1S6 CANADA
	\and 
	Evelyn Nelson\\Department of Mathematics and Statistics\\McMaster
	University\\Hamilton, Ontario L8S 4K1, CANADA \and Saharon
	Shelah\\Institute of Mathematics\\The Hebrew University\\Jerusalem, Israel}
	\date{}
	
\begin{document}
	\maketitle
	\begin{abstract}
	In the literature two notions of the word problem for a variety
	occur. A variety has a {\em decidable word problem} if every finitely
	presented algebra in the variety has a decidable word problem. It has
	a {\em uniformly decidable word problem} if there is an algorithm
	which given a finite presentation produces an algorithm for solving
	the word problem of the algebra so presented. A variety is given with
	finitely many axioms having a decidable, but not uniformly decidable,
	word problem. Other related examples are given as well.
	\end{abstract}

	{\S}0. INTRODUCTION 

	The following two options occur in the literature for what is meant by the 
	solvability of the word problem for a variety $V:$ 

	(1) There is an algorithm which, given a finite presentation ${\cal P}$ in 
	finitely many generators and relations, solves the word problem for
	${\cal P}$  relative to the variety $V$. 

	(2) For each finite presentation ${\cal P}$ in finitely many generators and 
	relations, there is an algorithm which solves the word problem for ${\cal P}$ 
	relative to the variety $V$. 

	We say that $V$ has uniformly solvable word problem if (1) holds.  It
	is the first notion that is studied in Evans [2, 3], where it is
	called just the word problem for $V$, and the second coincides with
	the terminology in Burris and Sankapannavar [1]. Benjamin Wells has
	informed us that Tarski was interested in the existence of varieties
	with solvable but not uniformly solvable word problem.  

	Varieties with uniformly solvable word problem include commutative
	semigroups and abelian groups (each of these are equivalent to the
	existence of an algorithm for solving systems of linear equations
	over the integers which is due to Aryabhata, see chapter 5 of [10]),
	any finitely based locally finite or residually finite variety, and
	the variety of all algebras of a given finite type (see [4]). 

	The examples which appear in the literature, of varieties with
	unsolvable word problem, all provide a finite presentation ${\cal P}$
	for which the word problem for ${\cal P}$ relative to that variety is
	unsolvable.  These include semigroups [9], groups [8], and modular
	lattices [5]. 

	Here, we present a finitely based variety $V$ of finite type which
	does not have a uniformly solvable word problem, but which
	nevertheless has solvable word problem.  We also present a
	recursively based variety of finite type, which is defined by laws
	involving only constants (i.e., no variables), with solvable but not
	uniformly solvable word problem.  This second result is the best
	possible one can provide for varieties defined by laws involving no
	variables: every finitely based such variety has uniformly solvable
	word problem.  If one would be satisfied with varieties with
	infinitely many operations, then it is relatively easy to produce an
	example of a recursively based variety with solvable but not
	uniformly solvable word problem; we present such an example,
	essentially due to B. Wells [11], at the end of the paper. 

	Our proof uses the unsolvability of the halting problem for the
	universal Turing machine, and the laws defining the variety precisely
	allow us to model the action of the universal Turing machine in the
	variety. In the usual proofs that the variety of semigroups has an
	undecidable word problem, a finitely presented congruence is given so
	that for any initial Turing machine configuration, the instantaneous
	descriptions of the Turing machine calculations all lie in the same
	congruence class. Then laws are added which make the halting state a
	right and left zero. So, the undecidability of the word problem in
	this algebra comes from not being able to decide whether a given word
	(initial configuration) is congruent to the halting state. This
	algebra contains all possible Turing machine calculations. In our
	variety each calculation will be modeled by a single algebra. 

	Our approach is based on a different picture of a Turing machine
	calculation than the sequence of instantaneous descriptions used in
	semigroups. We view a Turing machine calculation as taking place on a
	${\open Z} \times \omega $ grid, where the copy of ${\open Z}$ with
	second coordinate $n$ represents the Turing machine tape at time $n$.
	To understand the calculation, we must know the alphabet content of
	each square, which square the head is reading for each time $n$, and
	the state the machine is in at time $n$. There are various possible
	ways to formalize this insight, so that each Turing machine
	calculation corresponds to a finitely presented algebra. 

	To ensure that the word problem doesn't have a uniform solution, we
	introduce a function which has value 1 when applied to any state the
	machine reaches and value 0 on the halting state. Then a decision
	procedure which given a finite presentation, determines whether 1 is
	congruent to 0 would solve the halting problem. There are
	considerable technical difficulties in implementing this idea in such
	a way that we can prove that each finitely presented algebra in the
	variety has a decidable word problem.  

	It seems to us that there are two interesting directions that
	research can follow in light of the results in this paper. There
	remains the question of whether a finitely based variety of unary
	algebras with solvable word problem has uniformly solvable word
	problem. Another direction research could take is to consider
	subvarieties of interesting natural varieties. This could either be
	understood as varieties of $X$ where $X$ is a favourite class of
	algebras or, say, congruence modular varieties. This second problem
	was suggested to us by the persistent question of everyone to whom we
	told the result, namely ``Is there a natural example?" 

	The research for this paper was begun out while the latter two
	authors were visitors at the Department of Mathematics and Statistics
	at Simon Fraser University.  We gratefully acknowledge financial
	support from the Natural Science and Engineering Research Council of
	Canada.  This is paper \#291 on Shelah's publication list. We also
	wish to thank the referee for a thorough job of reading the paper. 

	\bigskip

	\noindent
	{\S}1. DEFINITIONS 

	\medskip

	We assume $\Sigma $ is some finitary type of algebras (with possibly
	infinitely many operations).  A presentation is a pair ${\cal P} =
	(X, R)$ consisting of a set $X$ (of generators) and a set $R
	\subseteq FX \times FX$ (of relations), where $FX$ is the (absolutely)
	free $\Sigma$-algebra over $X$.  A finite presentation is a
	presentation ${\cal P} = (X, R)$ where both $X$ and $R$ are finite.

	Given a variety $V$ and a presentation ${\cal P} = (X, R)$, there
	is an algebra $A \in V$ and a homomorphism $h: FX \rightarrow A$ with
	$R \subseteq Ker(h)$, such that any homomorphism $g: FX \rightarrow
	B$ with $B \in V$ and $R \subseteq Ker(g)$ factors uniquely through
	$h$.  The algebra A is unique up to isomorphism, and is called the
	algebra given by the presentation ${\cal P}$ relative to the variety
	$V$. 

	The word problem for ${\cal P}$ relative to the variety $V$ is to
	determine, given $s$, $t \in FX$, whether $(s, t) \in Ker(h)$.
	Note that $Ker(h)$ is the congruence on $FX$ generated by $R \cup
	\theta _V$, where $\theta _V$ consists of all equations in variables
	from $X$ satisfied by the variety $V$; equivalently, $\theta _V$ is
	the kernel of the unique homomorphism from $FX$ to the $V$-free algebra
	on $X$ mapping the elements of $X$ identically. 

	\medskip

	Next, we introduce the notion of a partial subalgebra, and state one result 
	which will be proved and used in {\S}6. The proof bears a familial
	resemblance to the more complicated proof in {\S}5.

	\medskip

	\noindent
	Definition.  A {\em partial subalgebra} is a pair $(A, \equiv _A)$ where $A$
	is a subset of $FX$ which is closed under formation of subterms, and
	$\equiv _A$ is an equivalence relation on $A$ which is a {\em partial
	congruence}, i.e., has the property that for each operation $\sigma $
	of arity $n$, if $a_i \equiv _A b_i$ for $1 \leq i \leq n$ and
	$\sigma (a_1,\ldots, a_n)$, $\sigma (b_1,\ldots, b_n) \in A$ then 
	$$\sigma (a_1,\ldots , a_n) \equiv _A \sigma (b_1,\ldots,
	b_n).  $$

	\medskip

	\noindent
	Proposition 1.1  {\it If for a partial subalgebra} $(A, \equiv _A)$,

	(1) {\it membership in} $A$ {\it is decidable (for elements of}
	$FX${\it ),} 

	(2) {\it membership in $\equiv _A$ is decidable (for pairs of elements of} 
	$FX${\it ),} 

	(3) {\it there is an algorithm which, given an operation $\sigma $ of
	arity $n$ and $a_1,\ldots$}, $a_n \in A$, {\it determines whether there
	exist $b_1,\ldots, b_n \in A$  for which $a_i \equiv _A b_i$ for
	$1 \leq i \leq n$ and $\sigma (b_1,\ldots, b_n) \in A$,}
	{\it then $\equiv ,$ the congruence on} $FX$ {\it generated by
	$\equiv _A$, is  
	decidable. Further this decision procedure is uniform in the algorithms for 
	deciding} (1), (2), (3){\it .}

	\medskip

	\noindent
	Remark: The conclusion of this result says that there is a solution
	to the word problem for the presentation $(X, \equiv _A)$ relative
	to the variety of all algebras of the given type.

	\medskip

	\noindent
	Definition A partial subalgebra satifying the hypotheses of
	Proposition 1.1 is called {\it decidable}. 

	\medskip

	We delay the proof of Proposition 1.1 until the end of {\S}6.

	\noindent
	Corollary 1.2. (Evans) {\it Let $V$ be the variety of all algebras in some 
	finite language.  Then $V$ has uniformly decidable word problem.} 

	\medskip

	Proof.  Suppose we are given some finite presentation.  Let $A$ be the
	finite set consisting of the terms appearing in the presentation and
	their subterms.  By brute search through the finitely many
	possibilities, we can find $\equiv _A$, the smallest partial
	congruence on $A$ containing all the relations in the presentation.
	Now we can apply Proposition 1.1. 

	\medskip

	Corollary 1.2 implies that any variety in a finite language which is
	defined by finitely many laws involving only constants has a
	uniformly decidable word problem. Any presentation can be viewed as a
	new presentation in the variety without the laws, by viewing each law
	as a relation in the new presentation. 

	\bigskip

	\noindent
	{\S}2. THE FINITELY BASED VARIETY.

	\medskip

	\noindent
	2.1 MODIFICATION OF THE UNIVERSAL TURING MACHINE 

	\medskip

	Suppose that we are given a universal Turing machine with a unique
	halting state $h$ which, at each
	move, prints some letter on the scanned square, moves one square
	either left or right (denoted respectively by $-1$ or $1$) and enters a
	new (or the same) state.  We are first going to adjust the machine by
	adding right and left end markers $e_R$ and $e_L$, (as new members of
	the alphabet), and adding, for each state, two new states $q_L$ and
	$q_R$, and appropriate instructions so that the adjusted machine does
	the following: if it is scanning $e_R$ in state $q$, it prints $B$
	(blank), moves right (into state $q_R)$, prints $e_R$ and then moves
	left and returns to state $q$; if it is scanning $e_L$ in state $q$
	it prints $B$, moves left, prints $e_L$, and moves right and returns
	to state $q$.  That is, to the Turing flow chart we add the following

	$$ \raisebox{2.5ex}{$q_L$}
	\stackrel{\textstyle\stackrel{e_L:B:-1}{\longleftarrow}}
	{\textstyle\stackrel{\textstyle\longrightarrow}{\scriptstyle e_L:B:1}}
	\raisebox{2.5ex}{$q$}
	\stackrel{\textstyle\stackrel{e_R:B:1}{\longrightarrow}}
	{\textstyle\stackrel{\textstyle\longleftarrow}{\scriptstyle e_R:B:-1}}
	\raisebox{2.5ex}{$q_R$}$$

	The resulting machine, if started on a finite tape inscription  with
	the left and right endmarkers at the appropriate ends, and the rest of
	the tape blank, does what the original machine would have done if
	placed on that inscription with the rest of the tape blank, except
	that whenever the adjusted machine hits an endmarker it first moves it
	out one square, leaving behind a blank square. 

	Suppose that the resulting machine has state set $Q$ and alphabet
	$\Sigma $, and that its action is given by the functions $\sigma $,
	$\mu $, and $\alpha $ operating on $Q \times  \Sigma $, which specify
	the next state, the motion (either left or right) and the print
	instruction, so that   

	$\sigma : Q \times  \Sigma  \rightarrow  Q$   

	$\mu : Q \times  \Sigma  \rightarrow  \{-1$, $1\}$   

	$\alpha : Q \times  \Sigma  \rightarrow  \Sigma .$

	For simplicity we will assume $\sigma$, $\mu$ and $\alpha$ are total
	functions and so defined even if we reach the ``halting state''.

	\medskip

	\noindent
	2.2 DEFINITION OF THE VARIETY

	\medskip

	\noindent
	Our variety $V$ has the following operations: 

	CONSTANTS: $c$, all elements of $Q \cup  \Sigma $, 0, 1, $0_F$, $1_F$ 

	UNARY:  $T$, $S$, $S^{-1}$, $H$, $P$, $C_\Sigma $, $C_Q$, $U$, $E$ 

	BINARY:  $F$, $R$, $K$, $K^*$, $C^*$,  

	TERNARY:  $N_H$, $N_Q$, $N_\Sigma $ 

	\medskip

	The elements of the range of $P$ will be the space-time elements (in
	the intended interpretation they represent the tape squares); the
	action of $S$ and $S^{-1}$ represents stepping right and left
	respectively though space, of $T$ represents moving ahead one time
	period, and of $H$ represents moving to the head position. $C_\Sigma$
	gives the letter in the square and $C_Q$ gives the state the machine
	is in while scanning the square. The action
	of the ternary operation $N_H$ gives head position at the next
	time instant, while $N_\Sigma$ gives the letter in the square scanned
	by the head at the next time position and $C_Q$ gives the state at
	the next time.  

	The intended interpretations of $F(x,y)$ and $R(x,y)$ are ``$y$
	follows $x$'' and ``$y$ is to the right of $x$''.  $K$  and  $K^*$ are
	comparison functions, and $U$ is a modified addition by 1. The
	intended interpretation is explained more fully in paranthetical
	comments below and in the proof of Theorem 3.1. 

	\medskip

	The laws defining our variety are the following:

	\noindent
	I \quad $PP(x) \approx  P(x)$  

	$PT(x) \approx  TP(x) \approx  T(x)$  

	$PS(x) \approx  SP(x) \approx  S(x)$  

	$PS^{-1}(x) \approx  S^{-1}P(x) \approx  S^{-1}(x)$  

	$PH(x) \approx  HP(x) \approx  H(x)$  

	$PN_H(x,y,z) \approx  N_H(x,y$, $P(z)) \approx  N_H(x,y,z)$  

	$PK(x,y) \approx  K(x,P(y)) \approx  K(x,y)$  

	\medskip

	$HS(x) \approx  HS^{-1}(x) \approx  HH(x) \approx  H(x)$  

	$HTH(x) \approx  HT(x)$  

	$TS(x) \approx  ST(x)$  

	$TS^{-1}(x) \approx  S^{-1}T(x)$  

	$SS^{-1}(x) \approx  S^{-1}S(x) \approx  P(x)$  

	\medskip

	$N_H(x,y$, $H(z)) \approx  N_H(x,y,z)$  

	$HN_H(x,y,z) \approx  HT(z).$

	\medskip

	\noindent
	II \quad  $N_Q(q,a,H(x)) \approx  \sigma (q,a)$  for all $q \in  Q$, $a \in  
	\Sigma $  

	$N_\Sigma (q,a, H(x)) \approx  \alpha (q,a)$  for all $q \in  Q$, $a \in  
	\Sigma $  

	$N_H(q,a, H(x)) \approx  S^{\mu (q,a)} TH(x)$  for all $q \in  Q - Q_{LR}$  

	\medskip

	$N_Q(q_L, C_\Sigma H(x)$, $H(x)) \approx  q \approx  N_Q(q_R, 
	C_\Sigma H(x)$, $H(x))$  

	\medskip

	$N_H(q_L, C_\Sigma H(x)$, $H(x)) \approx  S^{-1}TH(x)$  

	$N_H(q_R, C_\Sigma H(x)$, $H(x)) \approx  STH(x)$  

	\medskip

	$N_\Sigma (q_L, C_\Sigma H(x)$, $H(x)) \approx  e_L$  

	$N_\Sigma (q_R, C_\Sigma H(x)$, $H(x)) \approx  e_R$

	\medskip

	\noindent
	III \quad  $C_\Sigma  TH(x) \approx  N_\Sigma (C_QH(x), C_\Sigma H(x), H(x))$  

	$C_Q TH(x) \approx  N_Q (C_QH(x), C_\Sigma H(x), H(x))$  

	$HT(x) \approx  N_H(C_QH(x), C_\Sigma H(x), H(x))$

	\medskip

	\noindent
	IV \quad  $C_QP(x) \approx  C_QH(x) \approx  C_Q(x)$  

	$C_\Sigma (x) \approx  C_\Sigma P(x)$  

	$C_\Sigma TP(x) \approx  C^*(P(x), R(P(x), H(x)))$  

	$C_\Sigma TP(x) \approx  C^*(P(x), R(H(x), P(x)))$  

	$C_\Sigma P(x) \approx  C^*(P(x), 1)$

	($C^*$ ensures the symbol in a square, which is either to the ``right''
	or to the ``left'' of the square being scanned, remains unchanged at
	the next time.)
	\medskip

	\noindent
	V  \quad $R(x,y) \approx  R(P(x), P(y))$  

	$F(x,y) \approx  F(Px, Py)$  

	$R(P(x), P(x)) \approx  0$  

	$R(P(x), SP(y)) \approx  UR(P(x), P(y))$  

	\medskip

	$F(P(x), P(y)) \approx  F(H(x), H(y))$  

	$F(P(x), P(x)) \approx  0_F$  

	$F(P(x), TP(y)) \approx  UF(P(x), P(y))$  

	\medskip

	$U(0) \approx  U(1) \approx  1$  

	$U(0_F) \approx  U(1_F) \approx  1_F$

	\medskip

	\noindent
	VI \quad  For all operations $f$ except $T$, $S$, $S^{-1}$, $H$, $P$,
	$K$, and $N_H$, 
	$Pf$ is constant with value $P(c)$.

	\medskip
	\noindent
	VII (i)\quad $K(0, P(x)) \approx  P(x)$  

	$K(1, P(x)) \approx  P(c)$  

	$K(0_F, P(x)) \approx  P(x)$  

	$K(1_F, P(x)) \approx  P(c)$ 

	($K$ ensures that if $0 = 1$ or $0_F = 1_F$ in an algebra then all
	the space time elements are identical or in other words that
	space-time is degenerate.)
	\medskip

	(ii) \quad $K^*(d,d) \approx  0$ for all constants $d$  

	$K^*(d,e) \approx  1$ for all constants $d,e$ with $d \neq  e.$  

	$K^*(P(x), d) \approx  1$  for all constants  $d \neq  c$  

	$K^*(C_\Sigma P(x), d) \approx  1$ for all constants $d \notin  \Sigma $  

	$K^*(C_QP(x), d) \approx  1$ for all constants $d \notin  Q$  

	$K^*(R(P(x), P(y)), d) \approx  1$ for all constants $d \neq  0, 1$ 

	$K^*(F(P(x), P(y)), d) \approx  1$ for all constants $d \neq  0_F, 1_F$  

	$K^*(t,t) \approx  0$ and $K^*(s,t) \approx  1$ for all $s \neq  t$
	where both $s,t$ belong to    

	$\{P(x_1)$, $C_\Sigma  P(x_2)$, $C_Q P(x_3)$, $R(P(x_4)$, $P(y_4))$, 
	$F(P(x_5)$,  $P(y_5))\}$

	($K^*$ ensures that space-time is degenerate if there is any
	undesired equalities between constants or if there is a nonempty
	intersection between the ranges of certain operations.)

	\medskip
	\noindent
	VIII \quad  $EC_Q P(x) \approx  1$  

	$E(h) \approx  0$

	($E$ ensures that space-time is degenerate if the halting state is reached.)

	\medskip

	\noindent
	2.3 NORMAL FORM FOR SPACE-TIME ELEMENTS 

	The terms which are in the image of the operations  $P$, $S$,
	$S^{-1}$, $T$, $H$, $N_H$ and $K$  are called space-time terms.  For
	each such term $t$, $P(t)$ is equivalent, modulo the laws of our
	variety, to  $t$.  For all terms $t$  in the images of the other
	operations, $P(t)$ is equivalent, modulo the laws of our variety, to
	$P(c)$.  Thus a term  $t$  is a space-time term if and only if $t$
	and  $P(t)$ are equivalent,  modulo the laws of our variety. 

	We are going to develop a normal-form representation for space-time terms. 

	First, for a space-time term  $t$  define 
	\begin{eqnarray*}
	\Lambda _t  &=  &  \{S^nT^m(t) |m,n \in  {\open Z}, m \geq  0\} \cup
	\{S^nT^m HT^k(t)| n,m,k \in  {\open Z}, m,k \geq  0\} \\
	    &\cup  & \{S^nT^m N_H(s,u, T^k(t)) | n,m,k \in  {\open Z}, m,k \geq  0, 
	s, u \mbox{ arbitrary terms}\}
	\end{eqnarray*}

	\medskip

	Now, define the set  $G$  of {\em generating space-time terms} as follows 

	(i)   $P(c) \in  G$ 

	(ii)   $P(x) \in  G$  for each variable  $x$ 

	(iii)  For each term $t$, and each $g \in  G$ and $\lambda  \in
	\Lambda _g$, the term $K(t,\lambda ) \in  G.$ 

	(iv)   $G$  is the smallest set of terms satisfying (i), (ii) and (iii). 

	Further, define  

	$\Lambda  = \cup  \Lambda _g (g \in  G).$ 

	The members of $\Lambda _g$ for $g \in  G$ are called the space-time
	terms in normal form with respect to  $g$ and the $\Lambda_g$ is
	called the space-time component of $g$, in particular $g \in \Lambda_g$. 

	For space-time terms in normal form with respect to  $g$, we define the   
	$g$-time prefix, $g$-time coordinate and $g$-space coordinate as follows:

	\begin{tabular}{llll}
	term in $\gL$& $g$-time prefix & $g$-time coordinate& $g$-space
	coordinate\\ \hline   

	$S^nT^m(g)$& $T^m(g)$ & $m$ & $n$\\ 

	$S^nT^mHT^k(g)$& $T^mHT^k(g)$& $m + k$& $n$\\ 

	$S^nT^mN_H(s,t,T^k(g))$& $T^mN_H(s,t, T^k(g))$& $m + k + 1$& $n$
	\end{tabular}

	\medskip

	\noindent 
	Proposition 2.1:  {\it There is an effective procedure
	which, given a space-time term $s$}, {\it produces a term $t \in
	\Lambda$  (i.e., in normal form ) such that the laws} I {\it entail $s
	\approx  t.$} 

	Proof. The procedure is described inductively on the complexity of
	terms.  To begin, of course, the normal form of $g \in G$ is $g.$ If
	$t = P(s)$ for some normal form space-time term $s \in \Lambda _g$
	then the normal form for $t$ is the same as that of $s.$ 

	If $t = H(s)$ or $T(s)$, for some normal form space-time term $s \in
	\Lambda _g$, then the normal form $t'$ for $t$ is given in the
	following table:

	\begin{tabular}{lll}
	$s$ & $H(s)$ & $T(s)$ \\ \hline

	$S^nT^m(g)$ &    $HT^m(g)$ &$S^nT^{m+1}(g)$ \\

	$S^nT^m HT^k(g)$  & $HT^{m+k}(g)$ & $S^nT^{m+1} HT^k(g)$ \\

	$S^nT^mN_H(t_1,t_2,T^k(g))$ & $HT^{m+k+1}(g)$& $S^nT^{m+1}N_H(t_1,t_2,T^k(g))$
	\end{tabular}

	If $t = S(s)$ or $S^{-1}(s)$ for some normal form space-time term $s$
	then the normal form for $t$ is obtained from $s$ by adding or
	subtracting $1$ respectively to the space component. 

	If $t = N_H(s_1,s_2,s)$ for a normal form space-time term $s$, then
	the normal form $t'$ for $t$ is given in the following table: 

	\begin{tabular}{ll}

	$s$ &        $t'$\\ \hline

	$S^nT^m(g)$  & $N_H(s_1,s_2$, $T^m (g))$ \\

	$S^nT^m HT^k(g)$  & $N_H(s_1,s_2$, $T^{m+k} (g))$ \\

	$S^nT^m N_H(t_1,t_2, T^k(g))$ &  $N_H(s_1,s_2, T^{m+k+1} (g))$ 
	\end{tabular}

	If $t = K(s_1,s)$ for a normal form space-time term $s$  then $t$ is
	in normal form.

	This completes the description of the procedure.

	\medskip

	\noindent
	Remark.  In the ensuing development, we will always deal only with
	space-time elements in normal form, and when we write $H(\lambda )$,
	$S(\lambda )$, etc.  for $\lambda  \in  \Lambda $, we will mean the
	normal form of $H(\lambda )$, etc. 

	\medskip

	\noindent
	Proposition 2.2:  {\it For terms $s,t \in  \Lambda _g$ with time
	coordinates $m,n$ respectively, if $m < n$ then the laws $V$ entail
	$F(s,t) \approx  1_F.$} 

	\medskip

	\noindent
	Proposition 2.3:  {\it For terms $s,t \in  \Lambda _g$ with the same time 
	prefix but different space coordinates, the laws $V$ entail either $R(s,t) 
	\approx  1$ or $R(t,s) \approx  1.$} 

	\bigskip

	\noindent
	{\S}3. NON-UNIFORM SOLVABILITY OF THE WORD PROBLEM
	\medskip

	This section is devoted to a proof of the following:

	\medskip

	\noindent
	Theorem 3.1.  $V$ {\it does not have uniformly solvable word problem.} 

	\medskip

	Proof.  For any initial tape configuration 
	$$
	\begin{array}{ccccccccc}\hline
	\ldots\vline & e_L\vline & a_0\vline &a_1 \vline &a_2 \vline & \ldots
	\vline& a_k\vline & e_R\vline &\ldots\\ \hline
	 &   & \uparrow &
	\end{array}$$

	\noindent
	where $\uparrow $ indicates head position) for the universal Turing machine, 
	there is a corresponding finite presentation  

	${\cal P}:  C_Q(c) \approx  q_0$ and $C_\Sigma (S^{-1}(c)) \approx
	e_L$ and 
	$C_\Sigma (c) \approx  a_0$ and $\ldots C_\Sigma (S^k(c)) \approx
	a_k$ and $C_\Sigma (S^{k+1}(c)) \approx  e_R$  

	We claim that the universal Turing machine, started on that
	configuration, eventually halts, if and only if $E(q_0) \approx  h$
	(equivalently, $0 \approx 1)$ follows from the presentation ${\cal P}$
	in the variety $V$.  Thus, since there is no algorithm which
	determines, given an initial tape configuration, whether or not the
	universal Turing machine will halt, this establishes the fact that $V$
	does not have uniformly solvable word problem. 

	( $\rightarrow  ):$  This direction is clear; ${\cal P}$ together with
	the equations defining $V$ entail the analogous information at each
	successive configuration.  If the machine halts at time $n$ then we
	obtain $\lambda  \in \Lambda _{P(c)}$ such that $C_Q(\lambda ) \approx
	h$  and so   

	$$0 \approx  E(h) \approx  E(C_Q(\lambda )) \approx  1 \approx  E(q_0)$$
	follow from $\cal P$ in the variety $V$.
	\medskip

	( $\leftarrow  ):$  Suppose the machine, started on the above
	configuration, never halts.  Then we produce a model $A \in  V$
	satisfying all the equations in ${\cal P}$, in which $0 \neq  1.$ 

	The set of elements of  $A$  is  $\{*\} \cup  \Sigma  \cup  Q \cup
	\{S^nT^m(c)| n \in  {\open Z}, m \in  {\open N}\} \cup  \{n | n \in
	{\open Z} \mbox{ and } n \leq  1\} \cup \{n_F | n \in  {\open Z}
	\mbox{ and } n \leq  1\}.$ 

	The operations are defined in $A$ as follows: 

	(i)  $T,S,S^{-1}$ are defined on elements of the form  $S^nT^m(c)$
	according to equations I so as to yield elements again of this form;
	for other elements $y$, $T(y) = T(c)$, $S(y) = S(c)$, $S^{-1}(y) =
	S^{-1}(c).$ 

	\smallskip

	(ii)  $P$ maps all elements of the form $S^nT^m(c)$ identically and
	all other elements to $c$, in particular $P(c) = c$.

	\smallskip 

	(iii) $U(n) = (n+1)$ and $U(n_F) = (n+1)_F$ for $n \leq  0$, $U(1) = 1$, 
	$U(1_F) = 1_F$.  $U$ maps all other elements to $*$. 

	\smallskip 

	(iv)  $E(q) = 1$ for all $q \in  Q$, $q \neq  h$   

	$E(h) = 0$   

	$E$ maps all other elements to $*$. 

	\smallskip 

	(v)  $R(S^nT^m(c)$, $S^kT^j(c)) = \left\{  \begin{array}{ll}
			1&\hbox{ if }\ \ k > n\\ 
			k - n & \hbox{    if } k \leq  n
				\end{array}
	\right. $.   

	$R(x,y) = R(P(x), P(y))$ otherwise. 

	\smallskip 

	(vi)  $F(S^nT^m(c), S^kT^j(c)) =\left\{ 
			\begin{array}{ll}
			1_F & \hbox{ if } j > m\\
			(j - m)_F & \hbox{ if } j \leq   m
			\end{array}
	\right. $   

	$F(x,y) = F(P(x), P(y))$ otherwise. 

	\smallskip 

	(vii) $K(0, S^nT^m(c)) = S^nT^m(c) = K(0_F,S^nT^m(c))$   

	$K(1,S^nT^m(c)) = c = K(1_F,S^nT^m(c))$   

	$K(x,y) = c$ otherwise. 

	\smallskip 

	(viii) $K^*(d,d) = 0$ for all $d \in  \{0,1,0_F,1_F,c\} \cup  Q \cup
	\Sigma $  

	$K^*(d,e) = 1$ for all $d,e$ as above with $d \neq  e$   

	$K^*(k,d) = 1$ for all $k \in  \{n|n\leq 1\}$, all $d \neq  0,1$   

	$K^*(k,d) = 1$ for all $k \in  \{n_F|n \leq 1\}$, all $d \neq  0_F$, $1_F$   

	$K^*(S^nT^m(c), d) = 1$ for all constants $d \neq  c$ (including all
	$n$ and $n_F$, with  $n \leq  1)$   

	$K^*(x, y) = *$ else. 

	The values of  $H(x)$, $C_\Sigma (x)$, $C_Q(x)$, 
	$N_\Sigma (q,a,H(x))$, $N_Q(q,a,H(x))$, and $N_H(q,a,H(x))$ for  $x \in  
	\{S^nT^m(c)|n \in  {\open Z}, m \in  {\open N}\}$ are defined by induction on
	$m$: 

	define $H(c) = HS^n(c) = c$   

	$C_Q(c) = C_QS^n(c) = q_0$   

	$C_\Sigma (S^n(c))$ as in ${\cal P}$ for $-1 \leq  n \leq  k + 1$   

	$C_\Sigma (S^n(c)) = B$ for all other values of $n$   

	$N_Q(q,a,c)$,  $N_H(q,a,c,)$ and $N_\Sigma (q,a,c)$ are defined as in
	equations II (note that $c = H(c)).$ 

	Suppose we have already defined, for all $n \in  {\open Z}$,   

	$H(S^nT^m(c)) = H(T^m(c)) = S^kT^mc$ for some $k$  

	$C_Q(S^nT^m(c)) = C_Q(T^m(c)) \in  Q$  

	$C_\Sigma (S^nT^m(c)) \in  \Sigma $ 

	and $N_H$, $N_Q$, $N_\Sigma $ for all triples $(q,a, HT^m(c))$, with 
	appropriate values, i.e., $im(N_\Sigma ) \subseteq  \Sigma $ etc.

	Then define for all $n \in  {\open Z}$  

	$H(S^nT^{m+1}(c)) = N_H(C_Q(HT^m(c))$, $C_\Sigma (TH^m(c)), HT^m(c))$  

	$C_Q(S^nT^{m+1}(c)) = N_Q(C_Q(HT^m(c))$, $C_\Sigma (HT^m(c)), HT^m(c))$  

	$C_\Sigma (THT^m(c)) = N_\Sigma (C_Q(HT^m(c)), C_\Sigma (HT^m(c)), HT^m(c))$  

	$C_\Sigma (S^nTHT^m(c)) = C_\Sigma (S^nHT^m(c))$ for all $n \neq  0$ 

	and then define $N_Q$, $N_\Sigma $, $N_H$ for all triples $(q,a,
	HT^{m+1}(c))$ according to the rules II.

	\smallskip 
	\noindent
	This completes the inductive definition. 

	Define  $N_H(x,y,z) = N_H(x,y,P(z))$ if the latter has already been defined. 

	Define  $N_H$ and $H$ on all other elements to have value $c.$ 

	Define  $C_Q(y) = C_Q(c) = q_0$ for all $y$ not of the form  $S^nT^m(c)$ 

	$C_\Sigma (y) = C_\Sigma (c) = a_0$ , for all
	$y$  not of the form  $S^nT^m(c)$ 

	$N_Q(x,y,z) = N_\Sigma (x,y,z) = *$ for all values not defined above. 

	Define $C^*(S^nT^m(c),n) =\left\{ \begin{array}{ll}
		 C_\Sigma (S^nT^{m+1}(c))& \hbox{ if } n \leq  0\\
		 C_\Sigma (S^nT^m(c)) & \hbox{ if } n = 1
	 				 \end{array}
	\right. $ ,
	$C^*(x,y) = *$ otherwise. 

	Then the resulting algebra $A$ satisfies all the laws of the variety and the 
	equations of the presentation ${\cal P}$, and $0 \neq  1$ in $A$.

	\bigskip

	\noindent
	{\S}4. SOLVABLE WORD PROBLEM IN THE DEGENERATE CASE 

	This section and the next are devoted to proving that $V$ has solvable word 
	problem. 

	Let ${\cal P}$ be a finite presentation on a generating set $X =
	\{x_1,\ldots,x_n\}$ and let $\theta _{\cal P}$ be the congruence on $FX$
	generated by the relations of ${\cal P}$ together with the
	substitution instances of the laws defining our variety $V$.  We must
	prove that $\theta _{\cal P}$ is decidable.

	\medskip

	\noindent Definition. ${\cal P}$  has {\it degenerate space-time} if
	$P(t) \theta _{\cal P} P(c)$  for all terms $t$. Note that, by the
	laws of $V$ if ${\cal P}$  has degenerate space-time then   $0 \theta
	_{\cal P} 1$  and  $0_F \theta _{\cal P} 1_F$. In fact, the laws of
	$V$  allow this conclusion to be drawn from any failure of the
	operations  $S$, $T$  to behave without loops.  Also the laws of  VII
	imply that if either   $0 \theta _{\cal P} 1$  or  $0_F \theta _{\cal
	P} 1_F$, then space-time is degenerate.  

	\medskip

	We first prove the the following.

	\medskip

	\noindent
	Theorem 4.1 {\it If  ${\cal P}$  has degenerate space-time then the word 
	problem for  ${\cal P}$  relative to our variety  $V$  is decidable.}

	\medskip

	\noindent
	Proof. In this case, in the presented algebra  $F(X)/\theta _{\cal
	P}$, all the operations  $P$, $T$, $S$, $S^{-1}$, $H$ and $N_H$ are
	constant with value $P(c)$.  Moreover, $R$ and $F$ are constant, with
	value  $0$  and  $0_F$ respectively, $C_Q$ is constant with value
	$C_QP(c)$ and $C_\Sigma $ is constant, with value $C_\Sigma P(c)$.  

	\medskip

	Now, consider the type obtained from the one with which we are
	working, by deleting the operations  $P$, $T$, $S$, $S^{-1}$, $H$,
	$N_H$, $R$, $F$, $C_Q$, $C_\Sigma $, and adding three constants $c_1$,
	$c_2$, and $c_3$, which will stand for $P(c)$, $C_QP(c)$, and
	$C_\Sigma P(c)$ respectively.  Then there is an effective procedure
	which, given a term of the larger type, produces a term of the smaller
	type which is equivalent to it modulo the laws for our variety, and
	the equations $0 \approx  1$, $0_F \approx  1_F$, $c_1 \approx  P(c)$,
	$c_2 \approx  C_QP(c)$, $c_3 = C_\Sigma P(c).$ 

	\medskip

	Thus if we consider the variety $V'$ of this reduced type defined by the 
	following equations: 

	equations II for $N_Q$ and $N_\Sigma $, with space-time terms, and
	terms in the image of $C_Q$ and $C_\Sigma $ replaced by $c_1$, $c_2$,
	$c_3$ respectively. 

	$c_3 \approx  C^*(c_1$, 1) 

	$U(0) \approx  U(1) \approx  1$ 

	$U(0_F) \approx  U(1_F) \approx  1_F$ 

	equations VII, where $P(x)$ is replaced by $c_1$, $C_QP(x)$ by $c_2$, 
	$C_\Sigma P(x)$ by $c_3$ and $R(P(x)$, $P(y))$ and $F(P(x)$, $P(y))$ by 0. 

	$E(c_2) \approx  1$ 

	$E(h) \approx  0.$ 

	Then, if the terms in the presentation ${\cal P}$ are replaced by
	their equivalents in the new type we obtain a presentation ${\cal
	P}'$ which relative to the variety described above, is equivalent to
	${\cal P}$ relative to the original variety. Here by ``equivalent''
	we mean that there is an effective translation between terms so that
	${\cal P}'$ entails $s'\approx t'$ relative to the variety $V'$ if
	and only if $\cal P$ entails $s \approx t$ relative to the variety
	$V$. Since ${\cal P}'$ is a presentation considered relative to a
	variety defined by finitely many laws which involve no variables, 
	by Corollary~1.2 this word problem is solvable, which shows that the word
	problem for ${\cal P}$ relative to our variety is solvable too.

	\bigskip

	\noindent
	{\S}5. SOLVABLE WORD PROBLEM IN THE NON-DEGENERATE CASE

	\medskip

	\noindent
	5.1 PLAN OF THE PROOF 

	\medskip

	In this section we prove the following:

	\medskip

	\noindent
	Theorem 5.1 {\it If  ${\cal P}$  has non-degenerate space-time then
	the word problem for  ${\cal P}$  relative to our variety  $V$  is
	decidable.}

	\medskip

	\noindent
	Proof: The proof is presented in the remaining subsections of this
	section. In the remainder of this subsection we will describe the
	strategy of the proof. We will define by induction an increasing
	sequence of partial subalgebras  $(A_n, \equiv _{A_n})$. There are
	three points to be verified. First, for all $n$ every instance of the
	laws of the variety and the relations of  ${\cal P}$ with elements of
	$A_n$ is validated by  $\equiv _n$. Second, for  $a, b \in A_n$, if
	$a \equiv _n b$  then the laws of the variety and  ${\cal P}$  imply
	$a \theta _{\cal P} b$. The set of terms will, apart from passing to
	normal forms, equal  $\cup  A_n$. Hence  $\theta _{\cal P}$ will
	essentially equal $\cup  \equiv _n$. Third, the construction of the
	$A_n$ and the  $\equiv _n$ is uniformly effective and hence  $\cup
	\equiv _n$ and  $\theta _{\cal P}$ are decidable. We will give a
	careful definition of the  $A_n$ and  $\equiv _n$, but we will leave it
	to the reader to verify the three points mentioned above. One other
	point which is worth mentioning is that before constructing  $(A_0,
	\equiv _0)$  we will demand more information about  $\theta _{\cal
	P}$ other than its being non-degenerate. We first define two auxiliary
	sets  $A$  and  $B$.

	\medskip

	\noindent
	5.2 DEFINITION OF  $A$ 

	\medskip

	Let  $B_{\cal P}$ consist of all terms appearing in the presentation
	${\cal P}$ and all their subterms and all constants of the variety
	$V$.  Let $G_{\cal P}$ consist of  $P(c)$, $P(x)$ for each generator
	$x$ of the presentation, and all terms of the form  $K(s,t) \in
	B_{\cal P}$; thus $G_{\cal P}$ is finite. 

	For each $g \in  G_{\cal P}$, we let $\Gamma _{g,{\cal P}}$ be the set
	of all terms built from $g$ using the unary operations
	$P,S,S^{-1},T,H$, and $N_H(s,t,-)$ where $s,t \in  B_{\cal P}$.
	Further let $\Lambda _{g,{\cal P}}$ be the set of members of $\Gamma
	_{g,{\cal P}}$ which are in normal form with respect to $g$, and let
	$\Lambda _{\cal P} = \cup \Lambda _{g,{\cal P}} (g \in  G_{\cal P}).$ 

	Now, we define $A$ as follows:  it contains

	\noindent
	(1) the terms in $G_{\cal P}$ and all constants

	\noindent
	(2) the terms in $\Lambda _{\cal P}$

	\noindent
	(3) $U(d)$ for $d = 0,1,0_F,1_F$

	\noindent
	(4) $R(\lambda ,\gamma )$, $F(\lambda ,\gamma )$, $UR(\lambda ,\gamma )$, 
	$UF(\lambda ,\gamma )$ for $\lambda $, $\gamma  \in  \Lambda _{\cal P}$

	\noindent
	(5) $N_Q(q,a$, $HT^n(g)),$ 

	$N_\Sigma (q,a$, $HT^n(g)),$ 

	$N_H(q,a$, $HT^n(g))$ for all $g \in  G_{\cal P}$, $q \in  Q$, $a \in
	 \Sigma $ 

	$N_Q(q_L$, $C_\Sigma  HT^n(g)$, $HT^n(g)),$ 

	$N_Q(q_R$, $C_\Sigma  HT^n(g)$, $HT^n(g))$ for all $g \in  G_{\cal
	P}$, all $q \in  Q - Q_{LR}$, all $n \geq  0.$

	\noindent
	(6) $C_\Sigma (\lambda )$, $C_Q(\gamma )$, all $\lambda $, $\gamma  \in  
	\Lambda _{\cal P}$

	\noindent
	(7) $N_\Sigma (C_Q HT^n(g)$, $C_\Sigma  HT^n(g)$, $HT^n(g)),$ 

	$N_Q(C_Q HT^n(g)$, $C_\Sigma  HT^n(g)$, $HT^n(g)),$ 

	$N_H(C_Q HT^n(g)$, $C_\Sigma  HT^n(g)$, $HT^n(g))$ for all $g \in
	G_{\cal P}$, all $n \geq  0.$

	\noindent
	(8) $C^*(\lambda, R(\lambda, H\lambda )),$ 

	$C^*(\lambda, R(H\lambda, \lambda )),$ 

	$C^*(\lambda ,1)$ for $\lambda  \in  \Lambda _{\cal P}$

	\noindent
	(9) $K^*(d,d)$ for constants $d$ 

	$K^*(d,e)$ for all constants $d,e$ with $d \neq  e$.   

	All substitution instances of terms in laws VII (i) and (ii) where
	$P(x)$ and $P(y)$ are replaced by arbitrary $\lambda $, $\gamma  \in
	\Lambda _{\cal P}.$

	\medskip

	\noindent
	(10) $EC_Q(\lambda )$, for all $\lambda  \in  \Lambda _{\cal P}$ 
	and $E(h)$ 

	\medskip

	Note that membership in $A$ is decidable.

	\medskip

	\noindent
	5.3 DEFINITION OF  $B$ 

	\medskip

	Before we can define  $B$, we need some preliminary results.

	\medskip

	\noindent
	Lemma 5.2  {\it For any  $\lambda  \in  \Lambda _{g,{\cal P}}$},
	$\{\gamma \in  \Lambda _{g,{\cal P}}|\gamma \ \theta _{\cal P} \lambda
	\}$  {\it is finite.}

	\medskip

	\noindent
	Proof.  If $\lambda ,\gamma  \in  \Lambda _{g,{\cal P}}$ and $\lambda
	\ \theta _{\cal P} \gamma $ then it follows from laws of $V$ and the
	non-degeneracy of ${\cal P}$ that $\lambda $ and $\gamma $ have the
	same time coordinate. (For example,  $S^nT^kg\ \theta _{\cal P}
	S^mT^iHT^rg$  implies $HT^kg = H(S^nT^kg) \theta _{\cal P}
	H(S^mT^iHT^rg) = HT^{i+r}g$  and this yields  $k = i +r.)$ Moreover,
	two terms in $\Lambda _{g,{\cal P}}$ with the same time prefix and
	different space coordinates cannot be congruent modulo $\theta _{\cal
	P}$. Since there are only finitely many time prefixes with the same
	time coordinate as $\lambda $, this establishes the result.

	\medskip

	\noindent
	Corollary 5.3  {\it For any term  $t$}, $\{\lambda  \in  
	\Lambda _{\cal P}|\lambda \ \theta _{\cal P} t\}$ {\it is finite.}

	\medskip

	\noindent
	Definition.  For a finite $F \subseteq \Lambda _{\cal P}$, the {\em
	maximum time vector} of $F$ is $(m_g)_{g \in G_{\cal P}}$ where $m_g$
	is the {\em maximum $g$-time coordinate} of elements of $F \cap
	\Lambda _g$ (or 0 if $F \cap \Lambda _g$ is empty). We also make an
	{\em ad hoc} definition and say a space-time term $s$ is a {\em right
	subterm} of a space-time $t$ by induction on the construction of $t$.
	If $t$ is $Hu$, $Su$, or $S^{-1}u$ for a space-time term $u$ then $s$
	is a right subterm of $t$ if it is either $t$ or a right subterm of
	$u$. If $t$ is $N_H(w,v, u)$ where $u$ is a space-time term then $s$
	is a right subterm of $t$ if it is either $t$ or a right subterm of $u$.

	\medskip

	\noindent
	Lemma 5.4  {\it For any finite subset $F \subseteq  \Lambda _{\cal P}$ with 
	maximum time vector $(m_g)_{g \in  G_{\cal P}}$ there is a finite  $\bar F 
	\subseteq  \Lambda _{\cal P}$ with the same maximum time vector, such that} 

	(i) $F \subseteq  \bar F$ 

	(ii) {\it if $\lambda  \in  \Lambda _{\cal P}$ is a right subterm of
	$\gamma  \in   \bar F$ then $\lambda  \in  \bar F$} 

	(iii) {\it if $\lambda  \in  \Lambda _{\cal P}$ and $\lambda \ \theta
	_{\cal P} \gamma $ for $\gamma  \in  \bar F$ then $\lambda  \in  \bar
	F$} 

	(iv) {\it if $\lambda  \in  \Lambda _{\cal P}$ and the normal form of 
	$T\lambda $ belongs to $\bar F$ then $\lambda  \in  \bar F.$}

	\medskip

	\noindent
	Proof. We may assume that for each  $g \in  G_{\cal P}$, $T^jg \in
	F$ for all  $j \leq  m_g.$ 

	Now, let $g \in  G_{\cal P}$ and consider the set $F_g \se
	\Lambda _{g,{\cal P}} \cap  F$ which consists of all elements of 
	$\Lambda _{g,{\cal P}} \cap  F$ whose time coordinate relative to $g$
	is $m_g$. Let $F^*_g \se
	\Lambda _{g,{\cal P}}$ consist of all those $\lambda  \in  
	\Lambda _{g,{\cal P}}$ for which there exists $\gamma  \in  F_g$ with 
	$\lambda  \theta _{\cal P} \gamma $.  By Lemma 5.2, $F^*_g$ is finite. 

	Let $k_1$ and $k_2$ be the maximum and minimum, respectively, of the space 
	coordinate of members of $F^*_g$.  Note that, since $T^{m_g}(g) \in  F$, we 
	have $k_2 \leq  0 \leq  k_1.$ 

	Let $F'_g$ consist of those members of $\Lambda _g$ which are the
	normal forms of all terms of the form $S^k(\lambda )$ where $-k_1 \leq
	k \leq  -k_2$, and $\lambda  \in  F^*_g$.  Then $F'_g$ is finite and
	contains $F_g$.  We will show

	(a)   $\gamma  \in  F'_g$, $\lambda  \in  \Lambda _g$, 
	$\gamma \ \theta _{\cal P} \lambda $  implies  $\lambda  \in  F'_g$ 

	(b)   $\gamma  \in  F'_g$, $\lambda  \in  \Lambda _g$ a right subterm
	of $\gamma $  
	with time coordinate $m_g$ relative to $g$  implies  $\lambda  \in  F'_g.$ 

	re (a):  Suppose $\gamma $ is the normal form of $S^k(\delta )$ where
	$-k_1 \leq  k \leq  -k_2$ and $\delta  \in  F^*_g$.  Then $S^{-k}
	\gamma \theta _{\cal P} \delta $ and hence $S^{-k} \lambda  \theta
	_{\cal P} \delta $ and so the normal form of $S^{-k}\lambda $ belongs
	to $F^*_g$.  Thus $\lambda $, which is the normal form of
	$S^kS^{-k}\lambda $, belongs to $F'_g.$

	re (b): Suppose $\gamma \in F'_g$; then $\gamma $ is the normal form
	of a term $S^k(\delta )$ where $-k_1 \leq k \leq -k_2$ and $\delta
	\in F^*_g$.  Let the space coordinate of $\delta $ be $n$ and the
	time prefix of $\delta $ be $\tau $; then $\gamma = S^{n+k}\tau$.
	Moreover, all terms of the form $S^i\tau $ for $n - k_1 \leq i \leq n
	- k_2$ belong to $F'_g$.  Since $n - k_1 \leq 0 \leq n - k_2$, it
	follows that if $i$ is any number between $n + k$ and 0, then $S^i
	\tau \in F'_g$.  Now, any normal form subterm of $\gamma $ with the
	same time component has the same time prefix and hence this shows
	that every right subterm of $\gamma $ with the same time component belongs
	to $F'_g.$ 

	Let $F' = \cup F'_g (g \in G_{\cal P})$; then $F'$ is finite, $F
	\subseteq F'$, and $F'$ satisfies (ii) and (iii) for any $\lambda \in
	\Lambda _g$ with $g$-time coordinate $m_g$.  Add to $F'$ each term
	$\lambda \in \Lambda_g$ with $g$-time coordinate $m_g - 1$ such that
	the normal form of $T\lambda $ belongs to $F_g$.  The result is still
	finite.  Now repeat the procedure for elements of $g$-time coordinate
	$m_g - 1$ , etc. to eventually obtain the desired set $\bar F$.  This
	completes the proof.

	\medskip

	\noindent
	Definition. Define $B$ as follows:  

	Recall that $B_{\cal P}$ consists of all terms appearing in the
	presentation ${\cal P}$ and all subterms thereof, and all constants of
	our variety $V.$ 

	Enlarge $B_{\cal P}$ as follows: 

	(i)  For each $b \in  B_{\cal P}$, if there exists $a \in  A$ with $a 
	\theta _{\cal P} b$, add one such $a$, and choose $a \in  \Lambda _{\cal P}$ 
	whenever possible. 

	(ii)  For each $g \in  G_{\cal P}$ let $m_g$ be the maximum time
	coordinate of all the elements of $\Lambda _g$ that we have so far, and
	add all $g$-time prefixes up to time $m_g.$ 

	(iii) Let $F$ consist of all elements of $\Lambda _{\cal P}$ that we
	have so far.  Add the set $\bar F \supseteq  F$ given in the above
	lemma. 

	(iv)  For all $\lambda, \gamma  \in  \bar F$, add $C_Q(\lambda )$, 
	$C_\Sigma (\lambda )$, $R(\lambda ,\gamma )$, $F(\lambda ,\gamma )$, 
	$U(R(\lambda ,\gamma ))$, $U(F(\lambda ,\gamma )).$ 

	The resulting set $B$ is finite, is closed under taking subterms, and for 
	$\lambda $, $\gamma  \in  \Lambda_{\cal P} $, if $\gamma  \in  B$ and 
	$\lambda \ \theta _{\cal P} \gamma $ then $\lambda  \in  B.$ 

	Let $\equiv _B = \theta _{\cal P}|B$; then $\equiv _B$ is finite and hence 
	decidable. Also $\equiv_B$ contains the relations of $\cal P$.

	\medskip

	\noindent
	5.4 DEFINITION OF  $A_0$ 

	\medskip

	Now, define $A_0 = A \cup  B$; then membership in $A_0$ is decidable.
	We are going to define a partial congruence  relation $\equiv _0$ on
	$A_0$ so that the 
	pair $(A_0, \equiv _0)$ is a partial subalgebra such that membership
	in  $\equiv _0$ as well as  $A_0$ is decidable. In fact, $\equiv _0$
	will be $\theta _{\cal P}$ restricted to $A_0$, but we will define
	$\equiv _0$ by induction on the complexity of terms and the size of
	the time coordinate for members of $\Lambda _{\cal P}.$ 

	For $a,b \in  B$, $a \equiv _0 b$ if and only if $a \equiv _B b.$ 

	For $a \in  A$, $b \in  B$, $a \equiv _0 b$ if and only if there
	exists $c \in A \cap  B$ with $a \equiv _0 c$ (as described below) and
	$c \equiv _0 b$, i.e., $c \equiv _B b$.  Since $B$ is finite, we
	decide whether $a \equiv _0 b$ by searching through all $c \in  A \cap
	B$ and checking the latter two conditions.  Thus it is enough to
	describe $\equiv _0$ between pairs of elements of $A$. 

	There are some members of $A$ that we can essentially ignore, because
	we know they must be in the relation $\equiv _0$ to other elements
	that we have to deal with anyway. 

	Thus, to begin, we decree: 

	$U(0) \equiv _0 U(1) \equiv _0 1$ 

	$U(0_F) \equiv _0 U(1_F) \equiv _0 1_F$

	\noindent
	and for all $\lambda $, $\gamma  \in  \Lambda _p,$ 

	$UR(\lambda ,\gamma ) \equiv _0 R(\lambda ,\gamma ')$ where $\gamma '$ is the 
	normal form of $S(\gamma )$ 

	$UF(\lambda ,\gamma ) \equiv _0 F(\lambda ,\gamma ')$ where $\gamma '$
	is the normal form of $T(\gamma )$

	\noindent
	and so we may ignore, for the purposes of defining $\equiv _0$ between
	elements of $A$, all those elements of $A$ which are in the range of $U.$ 

	Similarly, using the appropriate terms given in laws II  and III we
	may ignore the elements of $A$ in the range of $N_\Sigma $ or $N_Q$, by
	making them $\equiv _0$ congruent to elements in $Q$ or the range of
	$C_Q$, and $\Sigma $ or the range of $C_\Sigma $, respectively. 

	We dispose in the same way of the elements of $A$ that are in the range
	of $C^*$, $K$ or $K^*.$ 

	Thus we need only define $\equiv _0$ between pairs of elements of $A$
	that are either constants, $\Lambda $-elements, or in the range of
	the operations $C_\Sigma, C_Q$, $F$ and $R$.  Moreover, since ${\cal
	P}$ is non-degenerate, we know by laws VIII that the interpreted
	images of these latter four operations are disjoint from one another
	and from all the interpretations of $\Lambda $-elements in
	$F(X)/\theta _{\cal P}$.  In addition, all the constants are pairwise
	distinct in $F(X)/\theta _{\cal P}$, and $A$-elements in the range of
	$C_\Sigma $, $C_Q$, $F$, and $R$ can be $\theta _{\cal P}$-congruent
	to constants only if they belong to $\Sigma $, $Q$, $\{0,1\}$,
	$\{0_F,1_F\}$, respectively.

	\medskip

	\noindent
	5.5 DEFINITION OF  $\equiv _0$ FOR ELEMENTS WITH SMALL TIME COMPONENT 

	\medskip

	We begin by describing $\equiv _0$ for elements $\lambda $, $C_\Sigma
	(\lambda )$, $C_Q(\lambda )$, $R(\lambda ,\gamma )$ and $F(\lambda
	,\gamma )$ for $\lambda $, $\gamma \in \Lambda _{\cal P}$ with time
	coordinate less than or equal to the maximum occurring in $B$,
	relative to whatever space-time component $\lambda $ and $\gamma $
	are in. 

	(1)  For $\lambda $, $\gamma  \in  \Lambda _{\cal P}$ with time coordinate 
	less than or equal to  the maximum in $B$, define   

	$\lambda  \equiv _0 \gamma $ if and only if $S^{-n} \lambda  \equiv
	_{\cal B} S^{-n} \gamma $

	\noindent
	where $n$ is the space coordinate of $\lambda .$ 

	Note that $S^{-n} \lambda  \in  B$, and hence if $\lambda \ \theta
	_{\cal P} \gamma $ then $S^{-n} \lambda \ \theta _{\cal P} S^{-n}
	\gamma $ and hence $S^{-n} \lambda  \equiv _B S^{-n} \gamma $.  The
	converse is also true, of course.  The point about the definition is
	that, given $\lambda $, we know its space coordinate and so we can
	decide $\lambda  \equiv _0 \gamma $ because $\equiv _B$ is decidable.
	Moreover (and we will need this later), given $\lambda $, we can
	calculate all (there are only finitely many) $\gamma  \in \Lambda
	_{\cal P}$ with $\lambda  \equiv _0 \gamma .$ 

	(2)  For $\lambda $, $\gamma  \in  \Lambda _{\cal P}$ with time coordinates 
	less than or equal to  the maximum in $B$,  

	(i) $R(\lambda ,\gamma ) \equiv _0 0$ if and only if
	\underline{either} $\lambda   
	\equiv _0 \gamma $ as in (1) above

	\underline{or} $\lambda, \gamma \in B$ and  $R(\lambda ,\gamma )
	\equiv _B 0.$   

	(ii) $R(\lambda ,\gamma ) \equiv _0 1$ if and only if 

	\underline{either} there
	exists $n > 0$ with $\gamma  \equiv_0  S^n\lambda $   

	\underline{or} there exists
	$\delta  \in  \Lambda _{\cal P} \cap  B$ and $n \geq 0$ with $\gamma
	\equiv _0 S^n \delta $ and $R(\lambda ,\delta ) \equiv _B
	0$  

	\underline{or} there exists $\delta  \in  \Lambda _{\cal P} \cap  B$
	and $n > 0$ with $\gamma  \equiv _0 S^n \delta $ and
	$R(\lambda,\delta ) \equiv _B 1.$  

	Note that this is decidable:  for example, to check whether there
	exists $n > 0$ with $\gamma  \equiv _0 S^n \lambda $, it is enough to
	determine whether there exists $n > 0$ with $S^{-m} \gamma  \equiv _B
	S^{n-m} \lambda $ where $m$ is the space coordinate of $\gamma $, and
	the latter is decidable because $B$ is finite.  

	(iii) $R(\lambda _1,\gamma _1) \equiv _0 R(\lambda _2,\gamma _2)$ if
	and only if

	\underline{either} both are congruent to 0 or 1 by (i) or (ii)   \ 

	\underline{or} $\lambda _1 \equiv _0 \lambda _2$ and $\gamma _1 \equiv _0 
	\gamma _2$   \ 

	\underline{or} there exist $\delta _1$ and $\delta _2$ and $n \geq  0$
	with $R(\lambda _1,\delta _1) \equiv _B R(\lambda _2,\delta _2)$ and
	$\gamma _1 \equiv _0 S^n \delta _1$, $\gamma _2 \equiv _0 S^n \delta
	_2.$

	\medskip

	\noindent
	Remark.  From the above definition, we have $R(\lambda ,\lambda )
	\equiv _0 0$ for all $\lambda  \in  \Lambda_{\cal P}$ with time
	coordinate less than 
	or equal to  the 
	maximum in $B$, and moreover, if $R(\lambda ,\gamma ) \equiv _0 0$ or
	1 then $R(\lambda ,S(\gamma )) \equiv _0 1$, and so the laws V for
	$R$ and these values of $P(x)$, $P(y)$, are satisfied. 

	(3)  For $\lambda $, $\gamma $, $\lambda '$, $\gamma ' \in  \Lambda
	_{\cal P}$ with time coordinate less than or equal to  the maximum in $B$, define  

	$F(\lambda ,\gamma ) \equiv _0 0_F$ if and only if $F(H\lambda $,
	$H\gamma ) \equiv _B 0_F$  

	$F(\lambda ,\gamma ) \equiv _0 1_F$ if and only if $F(H\gamma , H\delta ) 
	\equiv _B 1_F$  

	$F(\lambda ,\gamma ) \equiv _0 F(\lambda ',\gamma ')$ if and only if 
	$F(H\lambda, H\gamma ) \equiv _B F(H\lambda ', H\gamma ').$ 

	(4)  For $\lambda $, $\gamma  \in  \Lambda _{\cal P}$ with time coordinate 
	less than or equal to  the maximum in $B$, define  

	$C_Q(\lambda ) \equiv _0 q \in  Q$ if and only if $C_Q(H\lambda )
	\equiv _B q$  

	$C_Q(\lambda ) \equiv _0 C_Q(\gamma )$ if and only if
	\underline{either} both are  
	$\equiv _0$ the same $q \in  Q$    

	\underline{or} $C_Q(H\lambda ) \equiv _B C_Q(H\gamma ).$ 

	(5)  The description of when $C_\Sigma (\lambda ) \equiv _0 C_\Sigma
	(\gamma )$ is somewhat more complicated.  First, for space-time
	elements $\lambda $, $\gamma  \in  \Lambda $, define $\lambda
	\uparrow  \gamma $ to mean $\gamma  = T\lambda $ and either $R(\lambda
	,H\lambda ) \equiv _0 1$ or $R(H\lambda ,\lambda ) \equiv _0 1$.
	Further, define $\lambda  \downarrow \gamma $ to mean $\gamma
	\uparrow  \lambda .$ 

	Note that if $\lambda  \downarrow  \gamma  \uparrow  \delta $ then $\lambda  = 
	\delta $, and if $\lambda  \downarrow  \gamma  \equiv _0 \delta  \uparrow  
	\xi $ then $\lambda  \equiv _0 \xi .$ 

	Now define $\lambda  \uparrow ^* \gamma $ if and only if there is a
	finite sequence of $\uparrow $-moves from $\lambda $ to $\gamma $,
	i.e., if and only if there exist $\lambda _1$, $\lambda _2,\ldots,\lambda
	_k$ such that $\lambda  = \lambda _1 \uparrow  \lambda _2 \uparrow
	\lambda _3 \ldots\uparrow  \lambda _k = \gamma $.  Similarly define
	$\downarrow ^*.$ 

	Now, define  
	$$C_\Sigma (\lambda ) \equiv _0 C_\Sigma (\gamma )$$
	if and only if there exist natural numbers $k \leq  n$ and 
	$\lambda _1,\lambda _2,\ldots,\lambda _n \in \Lambda _p$ with time coordinates 
	less than or equal to  the maximum in $B$ such that 

	(i)  $\lambda  = \lambda _1$, or $C_\Sigma (\lambda ) \equiv _B 
	C_\Sigma (\lambda _1)$

	\noindent
	and (ii)  $\lambda _i \uparrow ^* \lambda _{i+1}$ for $i$ odd, $i \leq  k$   

	$\lambda _i \downarrow ^* \lambda _{i+1}$ for $i$ odd, $i > k$   

	$\lambda _i \equiv _0 \lambda _{i+1}$ for $i$ even

	\noindent
	and (iii) $\lambda _n = \gamma $ or $C_\Sigma (\lambda _n) \equiv _B 
	C_\Sigma (\gamma ).$ 

	Note that if such a sequence exists then the length of the shortest
	possible such sequence (including the lengths of the sequences
	involved in the $\uparrow ^*$ and $\downarrow ^*$ parts) is bounded above by
	twice the sum of the maximum time coordinates of elements in $B$.
	Hence we can decide, given $\lambda $ and $\gamma $, whether such a
	sequence exists. 

	Further, define  
	$$C_\Sigma (\lambda ) \equiv _0 a \in  \Sigma $$
	if and only if there exists $\gamma  \in  B$ with $C_\Sigma (\lambda ) 
	\equiv _0C_\Sigma (\gamma )$ as above and $C_\Sigma (\gamma ) \equiv _B a$. 

	With this definition, the congruence $\equiv _0$ (up to the maximum
	time coordinate in $B$) satisfies the laws IV.  The identities implied
	by laws II and III for $C_\Sigma $ are also satisfied because they
	only involve $C_\Sigma H(\lambda )$, and all the $H(\lambda )$ belong
	to $B.$

	\medskip

	\noindent
	5.6 COMPLETION OF THE DEFINITION OF  $\equiv _0$ 

	\medskip

	Now, we complete the definition of $\equiv _0$ for elements of $A$ which
	are, or which involve, space-time elements with time coordinate larger
	than the maximum in $B$, by induction on the time coordinate. 

	Suppose we have described $\equiv _0$ as above for pairs $(\lambda
	,\gamma )$, $(R(\lambda ,\gamma )$, $R(\lambda ',\gamma '))$, etc.\ 
	whenever the $\Lambda _{\cal P}$-elements $\lambda ,\gamma $, etc.\  in
	the space-time coordinate of $g$ have time coordinate  less than or
	equal $k_g$.  The following describes $\equiv _0$ for those elements
	in the space-time coordinate of $g$ involving time coordinate $k_g + 1.$

	1. (i)  For $\lambda ,\gamma  \in  A$, in the same space-time component, say 
	that of $g$, with time coordinate $k_g + 1$, define  

	$\lambda  \equiv _0 \gamma $ if and only if

	\underline{either} $\lambda  = \gamma $   

	\underline{or} $\lambda  = T\lambda '$, $\gamma  = T\gamma '$ and $\lambda ' 
	\equiv _0 \gamma '$   

	\underline{or} $\lambda  = S^n HT^{k_g+1} (g)$ and    

	either there exist $q \in  Q$, $a \in  \Sigma $ with  
	$C_Q HT^{k_g} (g) \equiv _0 q$ and $C_\Sigma  HT^{k_g} (g) \equiv _0 a$ 

	and $\gamma  = S^n T\ \gamma '$ for $\gamma ' \equiv _0 S^{\mu (q,a)}
	HT^{k_g} (g)$    

	or there exists $q \in  Q$ with $C_Q HT^{k_g}(g) \equiv _0 q_L$ and   
	$\gamma  = S^n T\ \gamma '$ for $\gamma ' \equiv _0 S^{-1} HT^{k_g}(g)$    

	or there exists $q \in  Q$ with $C_Q HT^{k_g}(g) \equiv _0 q_R$    

	and $\gamma  = S^nT\gamma '$ for $\gamma ' \equiv _0 SHT^{k_g}(g)$    

	or $\gamma  = S^n N_H(s,t,T^{k_g}(g))$ and $s \equiv _0 C_Q 
	T^{k_g}(g)$ 

	and $t \equiv _0 C_\Sigma  T^{k_g}(g).$   

	\underline{or} $\lambda  = S^n N_H(t_1,t_2$, $T^{k_g}(g))$, $\gamma  = S^n 
	N_H(s_1,s_2, T^{k_g}(g))$    \ \ 

	and $t_1 \equiv _B s_1$, $t_2 \equiv _B s_2$   

	\underline{or} vice-versa (with $\lambda ,\gamma $ switched).  

	(ii)  For $\lambda ,\gamma $ in different space-time components, say
	$\lambda $ in the space-time component of $g$ and $\gamma $ in the
	space-time component of $y$, with time coordinates $k_g + 1$ and less
	than or equal
	$k_y + 1$ respectively, define $\lambda  \equiv _0 \gamma $ if and
	only if one of the following holds:  

	Case 1  $\lambda  = S^n T^{k_g+1}(g)$ and there exists $\gamma '$ in the 
	space-time component of $y$ with $\gamma ' \equiv _0 T^{k_g}(g)$ and $\gamma  
	\equiv _0 S^n T\ \gamma '$ by the preceeding description for ``the same 
	space-time component".  

	Case 2  $\lambda  = S^n T^{m+1} HT^k(g)$ and there exists $\gamma '$
	in the space-time component of $y$ with $\gamma ' \equiv _0 T^m
	HT^k(g)$ and $\gamma \equiv _0 S^n T\ \gamma '.$  

	Case 3  $\lambda  = S^n HT^{k_g+1} (g)$ and there exists $\gamma '$ in
	the space-time component of $y$ with $\gamma ' \equiv _0 T^{k_g} (g)$
	and $\gamma \equiv _0 S^n H\ T\ \gamma '.$  

	Case 4  $\lambda  = N_H(t_1,t_2$, $T^{k_g}(g))$ and there exists
	$\gamma '$ in the space-time component of $y$ with  $\gamma ' \equiv
	_0 T^k(g)$ and $\gamma \equiv _0 N_H(t_1,t_2,\gamma ').$  

	Case 5  $\lambda  = S^n T^{m+1} N_H(t_1,t_2,T^k(g))$ and there exists 
	$\gamma '$ in the space-time component of $y$ with $\gamma ' \equiv _0 T^m 
	N_H(t_1,t_2,T^k (g))$ and $\gamma  \equiv _0 S^n T \gamma '.$ 

	2 (i) Define $R(\lambda ,\gamma ) \equiv _0 0$ where $\lambda $ has
	time coordinate $k_g + 1$ relative to $g$ and $\gamma $ has time
	coordinate less than or equal to $k_y + 1$ relative to $y$ (or
	vice-versa) if and only if $\lambda \equiv _0
	\gamma $ as in 1.  

	(ii)  Define $R(\lambda ,\gamma ) \equiv _0 1$ for $\lambda ,\gamma $
	as in (i) if and only if there exists $n > 0$ and $\delta  \in
	\Lambda _{\cal P}$ with $\gamma \equiv _0 S^n \delta $ and $R(\lambda ,\delta
	) \equiv _0 0$ (i.e., $\lambda \equiv _0 \delta$). 

	(Note that we can find all possible values for $\delta $ and hence can decide 
	whether these conditions are satisfied.)  

	(iii) Define $R(\lambda ,\gamma ) \equiv _0 R(\lambda ',\gamma ')$
	(where all of $\lambda ,\gamma ,\lambda ',\gamma '$ have time
	coordinate less than or equal to the relevant $k_g + 1$ and one of
	them has that time coordinate) if and only if either both $R(\lambda
	,\gamma )$ and $R(\lambda ',\gamma ')$ are $\equiv _0 0$ or both are
	$\equiv _0 1$ by (i) or (ii) respectively, or $\lambda \equiv _0
	\lambda '$ and $\gamma \equiv _0 \gamma '.$ 

	3. For $\lambda ,\gamma $ as in 2 we define  

	$F(\lambda ,\gamma ) \equiv _0 0_F$ if and only if $H\lambda  \equiv _0 
	H\gamma .$  

	Further, we will define $F(\lambda ,\gamma ) \equiv _0 1_F$ if and
	only if $F(H\lambda $, $H\gamma ) \equiv _0 1_F$, so it is enough to
	consider $\lambda = HT^k (g)$, $\gamma  = HT^m (y).$  

	If $k \leq  k_g$ and $m = k_y + 1$ then define   

	$F(\lambda ,\gamma ) \equiv _0 1_F$ if and only if $F(T^k(g), T^{k_y}(y)) 
	\equiv _0 0_F$ or $1_F.$  

	If $k = k_g + 1$ and $m \leq  k_y + 1$ then define   

	$F(\lambda ,\gamma ) \equiv _0 1_F$ if and only if there exists $n$
	with $0 < n \leq  m$ and $HT^k(g) \equiv _0 HT^{m-n}(y).$  

	Finally, define $F(\lambda ,\gamma ) \equiv _0 F(\lambda ',\gamma ')$
	if and only if either both are $\equiv _0$ to $0_F$ or both are
	$\equiv_0$ to $1_F$
	by the above, or $\lambda  \equiv _0 \lambda '$ and $\gamma  \equiv _0
	\gamma '$, or both $\lambda $ and $\lambda '$ have time coordinate
	less than
	the induction step, and there exists $n \geq  0$ with $H\gamma  =
	HT^{n+p} (g)$, $H\gamma ' = HT^{n+m} (y)$ and $F(\lambda , HT^p (g))
	\equiv _0 F(\gamma ', HT^m (y)).$ 

	4. For $\lambda $ in the space-time component of $g$ with time
	coordinate $k_g + 1$, we define $C_Q(\lambda ) \equiv _0 q \in  Q$ if
	and only if $C_Q(H\lambda ) \equiv _0 q \in  Q$, and for $H\lambda $,
	which is just $HT^{k_g+1}(g)$, we define $C_Q(H\lambda ) \equiv _0 q
	\in  Q$ if and only if either $C_Q HT^{k_g} (g) \equiv _0 q_R$ or
	$q_L$, or there exists $q' \in  Q$, $a \in  \Sigma $ with $\sigma
	(q',a) = q$ and $C_Q HT^{k_g} (g) \equiv _0 q'$ and $C_\Sigma
	HT^{k_g} (g) \equiv _0 a$.  

	Then, define $C_Q(\lambda ) \equiv _0 C_Q(\gamma )$ if and only if either 
	$H\lambda  \equiv _0 H\gamma $ or both $C_Q(\lambda )$ and $C_Q(\gamma )$ are 
	$\equiv _0$ to some $q \in  Q$ by the preceding paragraph. 

	5. Now, for $\lambda $ in the space-time component of $g$ with time
	coordinate $k_g + 1$, we define  

	$C_\Sigma (\lambda ) \equiv  a \in  \Sigma $ 

	if and only if

	\underline{Case 1}  $\lambda  = THT^{k_g} (g)$ and   

	\underline{either} there exists $q \in  Q$, $b \in  \Sigma $ with
	$\alpha (q,b) = a$ and 
	$$C_\Sigma (HT^{k_g} (g)) \equiv _0 b \mbox{ and } C_Q(HT^{k_g} (g)) \equiv _0 q$$   

	\underline{or} $a = e_R$ and there exists $q \in  Q$ with $C_Q(H\
	T^{k_g} (g)) \equiv _0 q_R$   

	\underline{or} $a = e_L$ and there exists $q \in  Q$ with $C_Q(H\
	T^{k_g} (g)) \equiv _0 q_L$   

	\underline{Case 2} There exists $\gamma $ with $\lambda  = T\gamma $
	and either  
	$R(\gamma ,H\gamma ) \equiv _0 1$ or $R(H\gamma ,\gamma ) \equiv _0 1$ and in 
	addition $C_\Sigma (\gamma ) \equiv _0 a$

	\underline{Case 3} There exists $\gamma $ in another space-time
	component with $C_\Sigma (\gamma ) \equiv _0 a$ by the above two cases
	and $\lambda  \equiv _0 \gamma .$  

	Finally, we define $C_\Sigma (\lambda ) \equiv _0 C_\Sigma (\gamma )$
	if and only if either both are $\equiv _0$ some $a \in  \Sigma $ by
	the preceding definition, or $\lambda  \equiv _0 \gamma $, or, if
	$\lambda $ is in the space-time component of $g$ with time coordinate
	$k_g + 1$ and $\gamma $ is in the space-time component of $y$ with
	time coordinate less than or equal to $k_y$ and there exists $\delta $
	with $\lambda  = T\delta $, and either $R(\delta ,H\delta )$ or
	$R(H\delta ,\delta ) \equiv _0 1$ and in addition $C_\Sigma (\delta )
	\equiv _0 C_\Sigma (\gamma )$ (or vice-versa with the roles of
	$\lambda $ and $\delta $ reversed). 

	This completes the definition of $(A_0, \equiv _0).$ As the notation
	suggests $\equiv_0$ is an equivalence relation, in fact a partial
	congruence. In particular transitivity is taken care of in the
	inductive construction. Furthermore it should be noted that $\equiv_0$
	is decidable. 

	\medskip\noindent
	5.7 DEFINITION OF  $(A'_n, \equiv '_n)$ 

	\medskip

	Now, suppose we have defined, for each $m \leq n$, a decidable
	partial subalgebra $(A_m, \equiv _m)$ such that $A_m \subseteq
	A_{m+1}$, $\equiv _m = \equiv _{m+1}|A_m$, $\equiv_m$ is a partial
	congruence, the laws of our variety are contained in $\equiv _m$
	insofar as they apply to the elements of $A_m$, and in addition 

	(i) $A_m$ is closed under  $P$, $H$, $S$, $S^{-1}$, $T$ (modulo normal
	form) and for all  $a \in  A_m$ and any operation among  $P$,
	$C_\Sigma $, $C_Q$, $R$, $F$, if  $a$  is  $\equiv _n$-equivalent to
	an element in the image of the operation then it is  $\equiv
	_m$-equivalent to an element in the image of that operation.

	\noindent
	Further we construct by induction algorithms which  

	(ii)  given space-time elements,  $\lambda $, $\gamma  \in   A_m$, determine 
	whether there exists  $k$  with  $\lambda  \equiv _m S^k\gamma .$    

	(iii) given space-time elements $\lambda $, $\gamma  \in  A_m$, determine 
	whether there exists  $k > 0$  with  $\lambda  \equiv _m T^k\gamma .$   

	(iv) given a space-time element  $\lambda  \in  A_m$ determine whether
	it is $\equiv _m$ to some element of  $A_0$,  and if so, produces such
	an element $\lambda '.$ 

	(v) given an element  $a \in  A_m$ and an operation among  $P$,
	$C_\Sigma $, $C_Q$, $R$, $F$, determine if  $a$  is $\equiv _m$ to an
	element in the range of the operation. 

	The base case of the construction is $n =0$. Note that $(A_0, \equiv
	_0)$ has the first four of these five properties: (i) follows from
	the definition of $A_0$ and the fact that $P(c) \in A_0.$ 

	(ii) is seen as follows, for  $\lambda , \gamma  \in  A_0:$ we can
	effectively list as  $\lambda _1, \ldots, \lambda _n$,  the finitely
	many elements in  $A_0$ to which  $\lambda$ is  $\equiv _0$. For some
	$k$, $\lambda  \equiv _0 S^k\gamma $  if and only if for some  $i$,
	$\lambda _i$ and  $\gamma $  have the same time prefix. The latter
	condition can be effectively checked. 

	(iii) Let  $\lambda _1, \ldots,
	\lambda _n$ be as above.  Then there is  $k$  so that  $\lambda
	\equiv _0 T^k\gamma $  if and only if some  $\lambda _i$ is in the
	space-time component of  $\gamma $,  $k_0$ is the difference of the
	time coordinates of  $\lambda _i$ and  $\gamma $,  and $\lambda
	_i \equiv_0  T^{k_0}\gamma .$ 

	(iv)  is trivial for $A_0.$

	\medskip

	\noindent Property (v) is a bit trickier. Since an element  $a$  is
	$\equiv _n$-equivalent to an element in the range of  $P$  if and only
	if  $a \equiv _n P(a)$, only the other operations present any
	problems. Below we will state an inductive hypothesis on the
	equivalence  $\equiv _n$. The inductive hypothesis will have two uses.
	First it will allow us to verify property (v) and the second part of
	(i) by giving a complete description of which elements are  $\equiv
	_n$-equivalent to an element in the range of a non-space-time
	operation. Second the hypothesis will determine  $\equiv '_n$, the
	restriction of  $\equiv _{n+1}$ to  $A'_n$ (defined below). In the
	remarks after the inductive hypothesis for elements in the image of
	$R$, we will expand on these points. The reader will be able to
	observe the inductive hypotheses hold for the case $n =0$ and so (v)
	holds as well. 

	\medskip

	Let  $A'_n$  be  $A_n$ together with the image of all the elements of
	$A_n$ under the operations  $R$, $F$, $U$, $C_Q$, $N_Q$, $C_\Sigma $,
	$N_\Sigma $, $C^*$, $E$, $K^*.$ 

	\medskip

	Extend  $\equiv _n$ to a partial congruence  $\equiv '_n$  on  $A'_n$
	by considering each operation in turn. First extend it to  $A_n$
	together with the image of  $R$  by letting it be the unique symmetric
	relation extending $\equiv _n$ which satisfies the inductive
	hypothesis. Continuing, given an operation we define  $\equiv '_n$  on
	the image of that operation,  $A_n$ and the operations previously
	considered, by letting it be the unique symmetric relation which
	satisfies the inductive hypothesis and extends the restriction of
	$\equiv '_n$  previously defined.  Transitivity will be an easy consequence
	of the definition since we will always link to a previous $\equiv_n$.

	\medskip

	(1) Image of $R$ 

	Inductive Hypothesis   

	$R(s, t) \equiv _n 0$ if and only if \underline{either} for some
	$\lambda _1$, $\lambda _2 \in A_0$, $P(s) \equiv _n \lambda _1$,
	$P(t)$ $\equiv _n \lambda _2$ and $R(\lambda _1$, $\lambda _2) \equiv
	_0 0$ \underline{or} $P(s) \equiv _n P(t).$ 

	$R(s, t) \equiv _n 1$ if and only if \underline{either} for some
	$\lambda _1$, $\lambda _2 \in A_0 P(s) \equiv _n \lambda _1$, $P(t)
	\equiv _n$ $\lambda _2$ and $R(\lambda _1, \lambda _2) \equiv _0 1$
	\underline{or} there is some $k > 0$ such that $P(t) \equiv _n
	S^kP(s).$ 

	$R(s, t) \equiv _n R(s', t')$ if and only if \underline{either} for
	some $\lambda _1$, $\lambda _2$, $\lambda _3$, $\lambda _4 \in A_0$,
	$P(s) \equiv _n \lambda _1$, $P(t) \equiv _n \lambda _2$, $P(s')
	\equiv _n \lambda _3$, $P(t') \equiv _n \lambda _4$ and $R(\lambda
	_1$, $\lambda _2) \equiv _0 R(\lambda _3$, $\lambda _4)$
	\underline{or} $R(s$, $t) \equiv _n 0 \equiv _n R(s'$, $t')$ (as
	above) \underline{or} $R(s, t) \equiv _n 1 \equiv _n R(s', t')$ (as
	above) \underline{or} $P(s) \equiv _n P(s')$ and $P(t) \equiv _n
	P(t').$ 

	$R(u, v) \equiv _n U(s)$ if and only if \underline{either} there are
	$t$, $\lambda _1$, $\lambda _2 \in B$ so that $t \equiv _n s$, $P(u)
	\equiv _n \lambda _1$, $P(v) \equiv _n \lambda _2$ and $U(t) \equiv
	_0 R(\lambda _1$, $\lambda _2)$ \underline{or} there are $\lambda
	_1$, $\lambda _2 \in A_n$ so that $s \equiv _n R(\lambda _1, \lambda _2)$ and
	$R(\lambda _1, S\lambda _2) \equiv _n R(u, v).$ 

	$R(s$, $t) \equiv _n u$  if and only if \underline{either} one of the above cases
	holds \underline{or} there is   $v \in  B$  such that  $R(s, t) \equiv _n v$  as
	above, and  $v \equiv _n u.$ 

	Remark. Note first that every equivalence on the right hand side of
	the inductive hypothesis either concerns elements in  $A_n$ or is
	covered in earlier clauses. As to the decidability of the relation, at
	various points we need to know if there are  elements with a certain
	property. For example in the fourth clause we ask whether ``there are
	$\lambda _1$, $\lambda _2$ so that $s \equiv _n R(\lambda _1,
	\lambda _2)$  and  $R(\lambda _1$, $S\lambda _2) \equiv _n R(u,
	v)$". By property (v) we can tell if there are $\lambda '_1$,
	$\lambda '_2$  so that  $s \equiv _n R(\lambda '_1, \lambda '_2)$.
	Since  $\equiv _n$ is a partial congruence, if there are $\lambda _1$,
	$\lambda _2$ so that  $s \equiv _n R(\lambda _1, \lambda _2)$ and
	$R(\lambda _1, S\lambda _2) \equiv _n R(u, v)$  then for all
	$\lambda '_1$, $\lambda '_2$,  $s \equiv _n R(\lambda '_1, \lambda
	'_2)$ implies that  $R(\lambda _1', S\lambda _2') \equiv _n R(u,
	v)$. Hence we have an algorithm for answering the question. Similar
	comments apply throughout. 

	There remains property (v) to consider. It is enough in view of the
	inductive hypothesis to be able to decide when elements of the form
	$U(s)$  and members of  $B$  are in the range of  $R$. Since  $B$  is
	finite, we can assume we know the answer for elements of  $B$
	relative to  $A_0$. By property (i), we will then know the answer for
	all  $A_n$, if we can settle the case $n =0$. By the fourth clause,
	we can reduce the question to either 
	elements of  $B$  or elements of  $A_0$. In  $A_0$ all the elements of
	the form  $U(s)$  are either in  $B$,  $\equiv _0$-equivalent to an
	element in the range of  $R$ or  $\equiv _0$-equivalent to an element
	in the range of  $F$. Since no element in the range of  $R$  can be
	$\equiv _0$-equivalent to an element in the range of  $F$, we can
	decide whether a given element of  $A_0$ in the range of  $U$ is
	$\equiv _0$-equivalent to an element in the range of  $R$. Such
	considerations recur throughout. 

	\medskip

	(2) Image of $F$ 

	Inductive Hypothesis  

	$F(s, t) \equiv _n 0_F$ if and only if \underline{either} there are
	$\lambda _1$, $\lambda _2 \in A_0$ so that $H(s) \equiv _n \lambda
	_1$, $H(t) \equiv _n \lambda _2$ and $F(\lambda _1, \lambda _2)
	\equiv _0 0_F$ \underline{or}  $H(s) \equiv _n H(t).$  

	$F(s, t) \equiv _n 1_F$ if and only if \underline{either} there are
	$\lambda _1$, $\lambda _2 \in A_0$ so that $H(s) \equiv _n \lambda
	_1$, $H(t)
	\equiv _n \lambda _2$ and $F(\lambda _1, \lambda _2) \equiv _0 1_F$
	\underline{or} there is $k > 0$ $H(t) \equiv _n HT^k(s).$ 

	$F(s, t) \equiv _n F(s', t')$ if and only if \underline{either} there
	are $\lambda _1$, $\lambda _2$, $\lambda '_1$, $\lambda '_2 \in A_0$
	so that $H(s) \equiv _n \lambda _1$, $H(t) \equiv _n \lambda _2$,
	$H(s') \equiv _n \lambda '_1$, $H(t') \equiv _n \lambda '_2$ and
	$F(\lambda _1, \lambda _2) \equiv _0 F(\lambda '_1, \lambda '_2)$
	\underline{or} $F(s, t) \equiv _n 0_F \equiv _n F(s', t')$ (as above)
	\underline{or} $F(s, t) \equiv _n 1_F \equiv _n F(s', t')$ (as above)
	\underline{or} $H(s)
	\equiv _n H(s')$  and  $H(t) \equiv _n H(t').$  

	$F(u, v) \equiv _n U(s)$ if and only if \underline{either} there are
	$t$, $\lambda _1$, $\lambda _2 \in B$ so that $t \equiv _n$ $s$,
	$H(u) \equiv _n \lambda _1$, $H(v) \equiv _n \lambda _2$ and $U(t)
	\equiv _0 F(\lambda _1$, $\lambda _2)$ \underline{or} there are
	$\lambda _1$, $\lambda _2$ so that $s \equiv _n F(\lambda _1$,
	$\lambda _2)$ and $F(\lambda _1, T\lambda _2) \equiv _n F(u, v).$ 

	$F(s, t) \equiv _n u$ if and only if \underline{either} one of the
	above cases holds \underline{or} there is $v \in B$ such that $F(s,
	t) \equiv _n v$ as above, and $v \equiv _n u$. 

	\bigskip

	(3) Image of $U$ 

	Inductive Hypothesis  

	$U(s) \equiv _n 0$ if and only if \underline{either} for some $t \in
	B$, $s \equiv _n t$ and $U(t) \equiv _0 0$ \underline{or} there are
	$\lambda _1$, $\lambda _2 \in A_n$ so that $s \equiv _n R(\lambda _1,
	\lambda _2)$ and $R(\lambda _1, S\lambda _2) \equiv _n 0$ 

	$U(s) \equiv _n 1$ if and only if \underline{either} for some $t \in
	B$, $s \equiv _n t$ and $U(t) \equiv _0 1$ \underline{or} there are
	$\lambda _1, \lambda _2 \in A_n$ so that $s \equiv _n R(\lambda _1,
	\lambda _2)$ and $R(\lambda _1, S\lambda _2) \equiv_n 1$ 

	$U(s) \equiv _n 0_F$ if and only if \underline{either} for some $t
	\in B$, $s \equiv _n t$ and $U(t) \equiv _0 0_F$ \underline{or} there
	are $\lambda _1$, $\lambda _2 \in A_n$ so that $s \equiv _n F(\lambda
	_1, \lambda _2)$ and $F(\lambda _1, T\lambda _2) \equiv_n 0_F$ 

	$U(s) \equiv _n 1_F$ if and only if \underline{either} for some $t
	\in B$, $s \equiv _n t$ and $U(t) \equiv _0 1_F$ \underline{or} there
	are $\lambda _1$, $\lambda _2 \in A_n$ so that $s \equiv _n F(\lambda
	_1, \lambda _2)$ and $F(\lambda _1, T\lambda _2) \equiv _n 1_F$ 

	$U(s) \equiv _n R(u, v)$ if and only if \underline{either} there are
	$t$, $\lambda _1, \lambda _2 \in B$, so that $t \equiv _n$ $s$, $P(u)
	\equiv _n \lambda _1$, $P(v) \equiv _n \lambda _2$ and $U(t) \equiv
	_0 R(\lambda _1, \lambda _2)$ \underline{or} there are $\lambda _1$,
	$\lambda _2 \in A_n$ so that $s \equiv _n R(\lambda _1, \lambda _2)$
	and $R(\lambda _1, S\lambda _2) \equiv _n R(u, v).$ 

	$U(s) \equiv _n F(u, v)$ if and only if \underline{either} there are
	$t$, $\lambda _1$, $\lambda _2 \in B$, so that $t \equiv _n s$, $P(u)
	\equiv _n \lambda _1$, $P(v) \equiv _n \lambda _2$ and $U(t) \equiv
	_0 F(\lambda _1$, $\lambda _2)$ \underline{or} there are $\lambda _1,
	\lambda _2 \in A_n$ so that $s \equiv _n F(\lambda _1, \lambda _2)$
	and $F(\lambda _1$, $T\lambda _2) \equiv _n F(u, v).$ 

	$U(s) \equiv _n U(t)$ if and only if \underline{either} $U(s) \equiv
	_n 0 \equiv _n U(t)$ (as above) \underline{or} $U(s) \equiv _n 1$
	$\equiv _n U(t)$ \underline{or} $U(s) \equiv _n 0_F \equiv _n U(t)$
	\underline{or} $U(s) \equiv _n 1_F \equiv _n U(t)$ \underline{or} $s
	\equiv _n t$ \underline{or} there are $s'$, $t' \in B$ so that $s
	\equiv _n s'$, $t \equiv _n t'$ and $U(s') \equiv _B U(t').$ 

	$U(s) \equiv _n u$ if and only if \underline{either} one of the above
	cases holds \underline{or} there is $v \in B$ such that $U(s) \equiv
	_n v$ as above, and $v \equiv _n u.$ 

	\bigskip

	(4) Image of $C_Q$ 

	Inductive Hypothesis   

	$C_Q(s) \equiv _n q$  for  $q \in  Q$  if and only if there is
	$\lambda  \in A_0$ such that  $H(s) \equiv _n \lambda $   and
	$C_Q(\lambda ) \equiv _0 q.$  

	$C_Q(s) \equiv _n C_Q(t)$ if and only if \underline{either} there are
	$\lambda _1$, $\lambda _2 \in A_0$ so that $H(s) \equiv _n \lambda
	_1$, $H(t) \equiv _n \lambda _2$ and $C_Q(\lambda _1) \equiv _n
	C_Q(\lambda _2)$ \underline{or} $H(s) \equiv _n H(t)$ 

	$C_Q(s) \equiv _n N_Q(u$, $v$, $w)$ if and only if \underline{either}
	$H(w) \equiv _n w$, $HT(w) \equiv _n H(s)$, $C_Q(w) \equiv _n u$ and
	$C_\Sigma (w) \equiv _n v$ \underline{or} for some $q \in Q$, $C_Q(s)
	\equiv _n q \equiv _n N_Q(u, v, w)$ \underline{or} there are $\lambda
	\in A_0$, and $u_0, v_0, w_0 \in B$ so that $H(s) \equiv _n \lambda
	$, $u_0
	\equiv _n u$, $v_0 \equiv _n v$, $w_0 \equiv _n w$ and $C_Q(\lambda )
	\equiv _0 N_Q(u_0, v_0, w_0).$ 

	$C_Q(s) \equiv _n u$ if and only if \underline{either} one of the
	above cases holds \underline{or} there is $v \in B$ such that $C_Q(s)
	\equiv _n v$ as above and $v \equiv _n u.$

	\noindent
	Remark. In the third clause we have to decide whether  $C_Q(w) \equiv
	_n u$  and $C_\Sigma (w) \equiv _n v$. In order that  $\equiv '_n$  be
	well defined we need to know that if  $N_Q(u, v, w) \in  A_n$ then
	$C_\Sigma (\lambda ) \in  A_n$ where  $\lambda $  is the normal form
	of  $H(w)$. This can be verified by induction on $n.$ 

	\medskip

	(5) Image of  $N_Q$ 

	Inductive Hypothesis  

	$N_Q(s, t, u) \equiv _n q \in Q$ if and only if \underline{either}
	there exist $p \in Q$ and $a \in \Sigma $ such that $s \equiv _n p$,
	$t \equiv _n a$ , $H(u) \equiv _n u$ and $\sigma (p, a) = q$
	\underline{or} $s \equiv _n q_L$ \underline{or} $q_R$, $H(u) \equiv
	_n u$ and $t \equiv _n C_\Sigma (u)$ \underline{or} $u \equiv _n
	H(u)$, $s \equiv _n C_Q(u)$, $t \equiv _n C_\Sigma (u)$ and
	$C_Q(T(u)) \equiv _n q$ \underline{or} there are $s'$, $t'$, $u' \in
	B$ so that $s \equiv _n s'$, $t \equiv _n t'$, $u \equiv _n u'$
	$N_Q(s,t,u) \in B$ and
	$N_Q(s'$, $t'$, $u') \equiv _B q.$ 

	$N_Q(s, t, u) \equiv _n N_Q(s', t', u')$ if and only if
	\underline{either} for some $q \in Q$, $N_Q(s, t, u) \equiv _n q
	\equiv _n N_Q(s', t', u')$ \underline{or} $s \equiv _n s'$, $t \equiv
	_n t'$ and $u \equiv _n u'$ \underline{or} $u \equiv _n H(u)$, $u'
	\equiv _n H(u')$, $s \equiv _n C_Q(u)$, $s' \equiv _n C_Q(u')$, $t
	\equiv _n C_\Sigma (u)$, $t' \equiv _n C_\Sigma (u')$ and $C_QT(u)
	\equiv _n C_QT(u')$ \underline{or} there exist $b \equiv _B b'$
	such that $N_Q(s, t, u) \equiv _n b$ (as above) and $N_Q(s', t', u')
	\equiv _n b'.$ 

	$N_Q(u, v, w) \equiv _n C_Q(s)$ if and only if \underline{either}
	$H(w) \equiv _n w$, $C_Q(w) \equiv _n u$, $C_\Sigma (w) \equiv _n v$
	and $C_QT(w) \equiv _n C_Q(s)$ \underline{or} for some $q \in Q$,
	$C_Q(s) \equiv _n q \equiv _n N_Q(u, v, w)$ \underline{or} there is
	some $b \in B$ so that $N_Q(u, v, w) \equiv _n b$ (as in the
	paragraph above) and $C_Q(s) \equiv _n b.$ 

	$N_Q(s, t, u) \equiv _n v$ if and only if \underline{either} one of
	the above cases holds \underline{or} there is $w \in B$ such that
	$N_Q(s, t, u) \equiv _n w$ as above, and $v \equiv _n w.$ 

	\medskip

	(6) Image of $C_\Sigma $ 

	Inductive Hypothesis 

	If  $R(s, H(s)) \equiv _n 1$  and  $P(s)$  is not equivalent to an
	element of $A_0$ then $R(T(s), HT(s))$  is not  $\equiv _n 1$  and
	dually for  $R(H(s), s).$  

	$C_\Sigma (s) \equiv _n a \in  \Sigma $  if and only if there are
	$\lambda \in  A_0$ such that  $\lambda  \equiv _n P(s)$  and
	$C_\Sigma (\lambda ) \equiv _0 a$.  

	$C_\Sigma (s) \equiv _n C_\Sigma (t)$ if and only if
	\underline{either} there are $\lambda _1$, $\lambda _2 \in A_0$ so
	that $P(s) \equiv _n \lambda _1$, $P(t) \equiv _n \lambda _2$ and
	$C_\Sigma (\lambda _1) \equiv _0 C_\Sigma (\lambda _2)$
	\underline{or} $P(s) \equiv _n P(t)$ \underline{or} $P(s) \equiv _n
	TP(t)$ and $R(t, H(t)) \equiv _n 1$ \underline{or} $R(H(t), t) \equiv
	_n 1$ \underline{or} the last condition holds with the roles of $s$
	and $t$ reversed.  

	$C_\Sigma (s) \equiv _n N_\Sigma (u, v, w)$ if and only if
	\underline{either} there is $a \in \Sigma $, so that $C_\Sigma (s)
	\equiv _n$ $a \equiv _n N_\Sigma (u, v, w)$ \underline{or} $u \equiv
	_n C_Q(w)$ and $v \equiv _n C_\Sigma (w)$, $w \equiv _n H(w)$ and
	$C_\Sigma (s) \equiv _n C_\Sigma (T(w))$ \underline{or} there are
	$\lambda \in A_0$, $u'$, $v'$, $w' \in B$ such that $P(s) \equiv _n
	\lambda $, $u \equiv _n u'$, $v \equiv _n v'$, $w \equiv _n w'$,
	$N_\Sigma (u', v', w') \in B$ and $C_\Sigma (\lambda ) \equiv _n
	N_\Sigma (u', v', w').$ 

	$C_\Sigma (s) \equiv _n C^*(t$, $u)$ if and only if
	\underline{either} for some $a \in \Sigma $, $C_\Sigma (s) \equiv _n
	a \equiv _n C^*(t, u)$ \underline{or} $u \equiv _n R(H(t), t)$
	\underline{or} $u \equiv _n R(t, H(t))$ and $C_\Sigma (s) \equiv _n
	C_\Sigma (T(t))$ \underline{or} $u \equiv _n 1$, $t \equiv _n P(t)$
	and $C_\Sigma (s) \equiv _n C_\Sigma (t)$ \underline{or} there are
	$t'$, $u' \in B$ and $\lambda \in A_0$ so that $t \equiv _n t'$, $u
	\equiv _n u'$, $C^*(t', u') \in B$, $P(s) \equiv _n
	\lambda $  and  $C^*(t', u') \equiv _0 C_\Sigma (\lambda ).$ 

	$C_\Sigma (s) \equiv _n u$ if and only if \underline{either} one of
	the above cases holds \underline{or} there is $v \in B$ such that
	$C_\Sigma (s) \equiv _n v$  as above and $v \equiv _n u.$ 

	\medskip

	(7) Image of $N_\Sigma $ 

	Inductive Hypothesis  

	$N_\Sigma (s, t, u) \equiv _n a \in \Sigma $ if and only if
	\underline{either} there are $q \in Q$, $b \in \Sigma $, such that $u
	\equiv _n H(u)$, $\alpha (q,b) = a$ and $s \equiv _n q$, $t \equiv
	_nb$ \underline{or} $a = e_L$ (resp.\ $e_R$) and there is $q \in Q$,
	such that $u \equiv _n H(u)$, $s \equiv _n q_L$ (resp. $q_R$), $t
	\equiv _n C_\Sigma (u)$ \underline{or} $H(u) \equiv _n u$, $C_Q(u)
	\equiv _n s$, $C_\Sigma (u) \equiv _n t$ and $a \equiv _n C_\Sigma
	T(u)$ \underline{or} there are $s'$, $t'$, $u' \in B$ so that $s
	\equiv _n s'$, $t \equiv _n t'$, $u \equiv _n u'$, $N_\Sigma (s', t',
	u') \in B$ and $N_\Sigma (s', t', u') \equiv _B a$.  

	$N_\Sigma (s, t, u) \equiv _n N_\Sigma (v, w, z)$ if and only if
	\underline{either} they are both $\equiv_n$-equivalent to the same $a
	\in \Sigma $ \underline{or} $s \equiv _n v$, $t \equiv _n w$ and $u
	\equiv _n z$ \underline{or} $H(u) \equiv _n u$, $C_Q(u) \equiv _n s$,
	$C_\Sigma (u) \equiv _n t$, $H(z) \equiv _n z$, $C_Q(z) \equiv _n v$,
	$C_\Sigma (z) \equiv _n w$ and $C_\Sigma T(u) \equiv _n C_\Sigma
	T(z)$ \underline{or} there are $b \equiv _B b' \in B$ so that
	$N_\Sigma (s, t, u) \equiv _n b$ and $N_\Sigma (v, w, z) \equiv _n
	b'$ (as above). 

	$N_\Sigma (u, v, w) \equiv _n C_\Sigma (s)$ if and only if
	\underline{either} there is $a \in \Sigma $, so that $C_\Sigma (s)
	\equiv _n$ $a \equiv _n N_\Sigma (u, v, w)$ \underline{or} $u \equiv
	_n C_Q(w)$ and $v \equiv _n C_\Sigma (w)$, $w \equiv _n H(w)$ and
	$C_\Sigma (s) \equiv _n C_\Sigma (T(w))$ \underline{or} there is
	$\lambda \in A_0$, $u'$, $v'$, $w' \in B$ such that $P(s) \equiv _n
	\lambda $, $u \equiv _n u'$, $v \equiv _n v'$, $w \equiv _n w'$,
	$N_\Sigma (u', v', w') \in B$ and $C_\Sigma (\lambda ) \equiv _n
	N_\Sigma (u', v', w').$ 

	$N_\Sigma (s, t, u) \equiv _n C^*(v, w)$ if and only if
	\underline{either} there is $a \in \Sigma $ so that $C^*(v, w) \equiv
	_n a \equiv _n N_\Sigma (s, t, u)$ \underline{or} there is
	some $z$ so that $C^*(v, w) \equiv _n C_\Sigma (z) \equiv _n N_\Sigma
	(s, t, u)$ \underline{or} there is some $b \in B$ so that $N_\Sigma
	(s, t, u) \equiv _n b$ (as above) and $b \equiv _n C^*(v, w).$ 

	$N_\Sigma (s, t, u) \equiv _n v$ if and only if \underline{either}
	one of the above cases holds \underline{or} there is $w \in B$ such
	that $N_\Sigma (s, t, u) \equiv _n w$ as above and $v \equiv _n w$. 

	\medskip

	(8)  Image of $C^*$ 

	Inductive Hypothesis  

	$C^*(u, v) \equiv _n a \in \Sigma $ if and only if \underline{either}
	there is $u 
	\equiv _n P(u)$, $v \equiv _n R(u, H(u))$ or $R(H(u), u)$ and
	$C_\Sigma T(u) \equiv _n a$ \underline{or} $v \equiv _n 1$, $u \equiv _n P(u)$
	and $P(u) \equiv _n a$ \underline{or} there exist $u'$, $v' \in B$ so that $u
	\equiv _n u'$, $v \equiv _n v'$, $C^*(u', v') \in B$ and  $C^*(u',
	v') \equiv _B a$.   

	$C^*(t, u) \equiv _n C_\Sigma (s)$ if and only if \underline{either}
	for some $a \in \Sigma$, $C_\Sigma (s) \equiv _n a \equiv _n C^*(t,
	u)$ \underline{or} $u \equiv _n R(H(t), t)$ \underline{or} $u \equiv
	_n R(t, H(t))$ and $C_\Sigma (s) \equiv _n C_\Sigma (T(t))$
	\underline{or} $u \equiv _n 1$, $t \equiv _n P(t)$ and $C_\Sigma (s)
	\equiv _n C_\Sigma (t)$ \underline{or} there are $t'$, $u' \in B$ and
	$\lambda \in A_0$ so that $t \equiv _n t'$, $u \equiv _n u'$,
	$C^*(t', u') \in B$, $P(s) \equiv _n \lambda $ and $C^*(t', u')
	\equiv _0 C_\Sigma (\lambda ).$ 

	$C^*(s, t) \equiv _n C^*(u, v)$ if and only if \underline{either}
	both are $\equiv_n$-equivalent to the same $a \in \Sigma$
	\underline{or} there are $\lambda _1$, $\lambda _2$ so that $C^*(s,
	t) \equiv _n C_\Sigma (\lambda _1)$, $C^*(u, v) \equiv _n C_\Sigma
	(\lambda _2)$ and $C_\Sigma (\lambda _1) \equiv _n C_\Sigma (\lambda
	_2)$ \underline{or} $s \equiv _n u$ and $t \equiv _n v$
	\underline{or} there are $b_1 \equiv _B b_2 \in B$ so that $C^*(s, t)
	\equiv _n b_1$ and $C^*(s, t) \equiv _n b_2$ as above.  

	$C^*(v, w) \equiv _n N_\Sigma (s, t, u)$ if and only if
	\underline{either} there is $a \in \Sigma $ so that $C^*(v, w) \equiv
	_n a \equiv _n N_\Sigma (s, t, u)$ \underline{or} there is
	some $z$ so that $C^*(v, w) \equiv _n C_\Sigma (z) \equiv _n N_\Sigma
	(s, t, u)$ \underline{or} there is some $b \in B$ so that $N_\Sigma
	(s, t, u) \equiv _n b$ and $b \equiv _n C^*(v, w)$ (as above).  

	$C^*(s, t) \equiv _n u$ if and only if \underline{either} one of the
	above cases holds \underline{or} there is $v \in B$ such that $C^*(s,
	t) \equiv _n v$ as above and $u \equiv _n v.$ 

	\medskip

	(9) Image of $E$ 

	Inductive Hypothesis  

	$E(s) \equiv _n 1$ if and only if \underline{either} for some
	$\lambda $, $s \equiv _n C_Q(\lambda )$ \underline{or} there is $t
	\in B$ so that $s \equiv _n t$ and $E(t) 
	\equiv _B 1.$  

	$E(s) \equiv _n 0$ if and only if \underline{either} $s \equiv _n h$
	\underline{or} there is $t 
	\in  B$  so that  $s \equiv _n t$   and  $E(t) \equiv _B 0.$  

	$E(s) \equiv _n E(t)$ if and only if \underline{either} $s \equiv _n
	t$ \underline{or} $E(s) \equiv _n 1 \equiv _n E(t)$ \underline{or}
	$E(s) \equiv _n$ $0 \equiv _n E(t)$ \underline{or} there are $b_1
	\equiv _B b_2$ so that $E(s) \equiv _n b_1$ and $E(t)
	\equiv _n b_2$ (as above).  

	$E(s) \equiv _n t$ if and only if \underline{either} one of the above
	cases holds \underline{or} there is $u$ $\in B$ such that $E(s)
	\equiv _n u$ (as above) and $u \equiv _n t.$ 

	\medskip

	(10) Image of $K^*$ 

	Inductive Hypothesis  

	$K^*(s, t) \equiv _n 0$ if and only if \underline{either} $s \equiv _nt$ and
	either $s$ is $\equiv_n$-equivalent to a constant or $P(s) \equiv _n
	s$ or $s$ is equivalent to an element in the image of
	$C_\Sigma $, $C_Q$, $R$ or $F$ \underline{or} there are $b_1$, $b_2
	\in B$ so that $s \equiv _n b_1$, $t \equiv _n b_2$ and $K^*(b_1,
	b_2) \equiv _B 0.$ 

	$K^*(s, t) \equiv _n 1$ if and only if \underline{either} $P(s)
	\equiv _n s$ and there is a constant $d \neq c$ with $t \equiv _n d$
	\underline{or} $s$ is equivalent to an element in the range of
	$C_\Sigma $ and there is a constant $d \notin \Sigma $ so that $t
	\equiv _n d$ \underline{or} $s$ is equivalent to an element in the
	range of $C_Q$ and there is a constant $d \notin Q$ so that $t \equiv
	_n d$ \underline{or} $s$ is equivalent to an element in the range of
	$R$ and there is a constant $d \neq 0, 1$ so that $t \equiv _n d$
	\underline{or} $s$ is equivalent to an element in the range of $F$
	and there is a constant $d \neq 0_F$, $1_F$ so that $t \equiv _n d$
	\underline{or} $s$ is equivalent to an element in the range of one of
	the operations $P$, $C_\Sigma $, $C_Q$, $R$, $F$ and $t$ is
	equivalent to an element in the range of a different one of these
	operations \underline{or} there are $b_1$, $b_2 \in B$ so that $s
	\equiv _n b_1$, $t \equiv _n b_2$ and $K^*(b_1, b_2) \equiv _B 1.$ 

	$K^*(s, t) \equiv _n K^*(u$, $v)$ if and only if \underline{either}
	$K^*(s, t) \equiv _n 0 \equiv _n K^*(u, v)$ \underline{or} $K^*(s, t)
	\equiv _n 1 \equiv _n K^*(u, v)$ \underline{or} $s \equiv _n u$ and
	$t \equiv _n v$ \underline{or} there is $b_1 \equiv _B b_2$ so that
	$K^*(s, t) \equiv _n b_1$ and $K^*(u, v) \equiv _n b_2$ (as above).  

	$K^*(s, t) \equiv _n u$ if and only if \underline{either} one of the
	above cases holds \underline{or} there is $v \in B$ such that $K^*(s,
	t) \equiv _n v$ (as above) and $u \equiv _n v.$

	\medskip

	\noindent
	5.8 DEFINITION OF  $(A^{''}_n, \equiv^{''}_n)$

	\medskip

	\noindent
	Define $A^{''}_n$ and $\equiv^{''}_n$ as follows:

	\noindent
	for $n$ even: 

	$A^{''}_n
	= A'_n \cup  \{N_H(s,t,H\lambda )|s,t \in  A'_n, \lambda  \in  A_n 
	\mbox{ space-time}\}$

	\noindent
	for $n$ odd: 

	$A^{''}_n
	= A'_n \cup  \{K(s,\lambda )|s \in  A'_n, \lambda  \in  A_n \mbox{
	space-time}\}$  

	\medskip

	We will extend $\equiv '_n$ to a partial congruence $\equiv^{''}_n $
	on $A^{''}_n $.  First, for $n$ even, we define $N_H(s,t,H\lambda )
	\in A^{''}_n -A_n$ to be reducible (to $u$) if and only if one of the
	following holds: 

	(i) there exists $s',t',\lambda '$ with $u = N_H(s',t',H\lambda ')
	\in B$ and $s \equiv' _ns'$, $t \equiv' _nt'$ and $H\lambda \equiv _n
	H\lambda '.$ 

	(ii) $s \equiv'_n q \in Q$, and $t \equiv'_n a \in \Sigma $, and $u =
	S^{\mu (q,a)}TH(\lambda )$ 

	(iii) $s \equiv' _n q_L$ and $t \equiv '_n C_\Sigma H(\lambda )$ and
	$u = S^{-1}TH(\lambda )$ 

	(iv) $s \equiv' _n q_R$ and $t \equiv '_n C_\Sigma H(\lambda )$, and
	$u = STH(\lambda )$ 

	(v)  $s \equiv '_n C_QH(\lambda )$, $t \equiv '_n C_\Sigma H(\lambda )$ and 
	$u = HT\lambda .$ 

	\medskip

	Now define, for $n$ even, 

	$N_H(s,t,H\lambda ) \equiv^{''}_n N_H(s',t',H\lambda ')$ (for both
	$N_H(s,t,H\lambda )$ and $N_H(s',t',H\lambda ') \in A^{''}_n -A_n)$
	if and only if \underline{either} $s \equiv '_n s'$, $t \equiv
	'_n t'$, $H\lambda \equiv H\lambda $ \underline{or} both
	$N_H(s,t,H\lambda )$ and $N_H(s',t',H\lambda ')$ are reducible, to
	$u$, $u'$ respectively , and $u \equiv' _n u'$ 

	$N_H(s,t,H\lambda ) \equiv^{''}_n v\ \epsilon \ A_n$ if and only if
	\underline{either} $N_H(s,t, H\lambda )$ is reducible to $u$, and $u$
	$\equiv' _n v$, \underline{or} $v = N_H(s',t', H\lambda ')$ and $s
	\equiv' _n s'$, $t \equiv' _n t'$, $H\lambda \equiv' _n H\lambda '.$ 

	Similarly, for $n$ odd, we define $K(s,\lambda ) \in  A^{''}_n
	 - A_n$ to be reducible (to $u$) if and only if
	one of the following hold 

	(i)  there exist $s',\lambda '$ with $u = K(s',\lambda ') \in B$ and $s 
	\equiv '_n s$, $\lambda  \equiv' _n \lambda '$  

	(ii)  $s \equiv  '_n 0$ or $0_F$ and $u = \lambda $  

	(iii) $s \equiv '_n 1$ or $1_F$ and $u = P(c).$ 

	\medskip

	Now , for $n$ odd define $K(s,\lambda ) \equiv^{''}_n K(s',\lambda
	')$ (for both $K(s,\lambda )$ and $K(s',\lambda )$ and $K(s',\lambda
	')$ in $A^{''}_n - A_n$ if and only if \underline{either} $s \equiv
	'_n s'$ and $\lambda \equiv' _n \lambda '$ \underline{or} both
	$K(s,\lambda )$ and $K(s',\lambda ')$ are reducible to $u,u'$
	respectively, and $u \equiv' _n u'.$ 

	$K(s,\lambda ) \equiv '_n v \in  A_n$ if and only if $K(s,\lambda )$ is 
	reducible to $u$, and $u \equiv' _n v$
	\underline{or} $v = K(s',\lambda ')$ and both $K(s,\lambda )$ and $K(s'$,
	$\lambda ')$ are irreducible, and 
	$s \equiv _n s'$, $\lambda  \equiv' _n \lambda .$ 

	\medskip

	Finally, we define $A_{n+1}$ and $\equiv _{n+1}$ to extend $A^{''}_n$
	and $\equiv^{''}_n$ as follows: 

	For $n$ even:  $A_{n+1} = A'_n \cup  \{S^nT^m\lambda | \lambda  \in  A^{''}_n
	- A_n$, $n \in  {\open Z}, m \geq  0\}.$ 

	For $n$ odd:  
	\begin{eqnarray*}
	A_{n+1} = A^{''}_n&\cup&  \{S^nT^m \lambda |\lambda  \in  A^{''}_n
	- A_n, n \in  {\open Z}, m \geq  0\} \\                           
	&\cup & \{S^nT^m HT^k \lambda |\lambda  \in  A^{''}_n -A_n, n \in
	{\open Z}, m,k \geq  0\}. 
	\end{eqnarray*}

	\medskip

	\noindent
	5.9 DEFINITION OF  $(A_{n+1}, \equiv _{n+1})$  

	\medskip

	We extend $\equiv^{''}_n$
	to $\equiv _{n+1}$ on $A_{n+1}$ as follows: 

	\medskip

	For $n$ even: 

	$S^nT^m \lambda  \equiv _{n+1} S^{n'}T^{m'}
	$ if and only if $n = n'$, $m = m'$ and $\lambda  \equiv^{''}_n
	\lambda '$ and  

	$S^nT^m\lambda  \equiv _{n+1} v \in  A_n$ if and only if $\lambda $ is 
	reducible, to $u$, and $S^nT^mu \equiv_n v.$ 

	For $n$ odd:   

	$S^nT^m\lambda \equiv _{n+1} S^{n'}T^{m'}\lambda '$ if and only if $n
	= n'$, $m = m'$, and $\lambda \equiv _n \lambda '.$ 

	$S^nT^mHT^k\lambda \equiv _{n+1} S^{n'}T^{m'}HT^{k'}\lambda '$ if and
	only if $n = n'$. $m = m'$, $k = k'$, and $\lambda = \lambda '.$ 

	Finally, $S^nT^m\lambda \equiv _{n+1} v \in A_n$ if and only if
	$\lambda $ reduces to $u$ and $S^nT^mu \equiv _nv$, and 
	similarly with $S^nT^mHT^k\lambda .$

	\medskip

	\noindent
	This completes the definition of $(A_{n+1}, \equiv _{n+1})$.  It is
	straightforward to check that $(A_{n+1}$, $\equiv _{n+1})$ is a
	decidable partial congruence satisfying (i), (ii), (iii), (iv) above.
	Property (v) can be verified in the other cases as it was after the
	definition of the inductive hypothesis on $R$.  Also, as discussed in
	subsection 5.1, $\bigcup_{n\geq 0}\equiv_n$ is precisely $\theta
	_{\cal P}$ restricted to $ \bigcup_{n\geq 0}A_n$. The $A_n$ and
	$\equiv_n$ are uniformly decidable in $n$. Together this means that
	$\bigcup_{n \geq 
	0} A_n$ and $\bigcup_{n \geq 0} \equiv_n$ are decidable.  Since $
	\bigcup_{n\geq 0}A_n$ contains all terms, modulo (effectively)
	reducing space-time terms to their normal forms, this completes the
	proof.

	\medskip

	If we combine Theoerems 3.1, 4.1 and 5.1, we obtain the following theorem.

	\medskip

	\noindent
	Theorem 5.5  {\it There is a finitely based variety of finite type which has 
	solvable but not uniformly solvable word problem.}

	\bigskip

	\noindent
	{\S}6. A RECURSIVELY BASED VARIETY DEFINED BY LAWS INVOLVING NO VARIABLES.

	\medskip

	\noindent
	6.1 DEFINITION OF THE VARIETY 

	\medskip

	In this section we will describe a recursively based variety of finite type, 
	defined by laws which involve no variables, which has solvable but not 
	uniformly solvable word problem. 

	The variety is a modification of the finitely based variety defined in the 
	preceding sections. The use of infinitely many axioms allows us to
	use a simpler picture of space-time.  

	The operations are the same, except that  $P$, $C^*$,  and  $U$  are deleted 
	and $K$ and $K^*$ are identified.  Specifically, the operations are:  

	\medskip

	constants:  $c$, all $a \in  \Sigma $, all $q \in  Q$, 0, 1, $0_F$, $1_F$  

	unary:      $T$, $S$, $S^{-1}$, $H$, $C_\Sigma $, $C_Q$, $E$  

	binary:     $F$, $R$, $K$  

	ternary:    $N_H$, $N_Q$, $N_\Sigma .$ 

	\medskip

	Define, for each $k$, $H_k = HT^k(c)$, and let $\Lambda  = \{S^nT^m(H_k) |n 
	\in  {\open Z}, m, k \in  {\open N}\}$.  These will be the space-time 
	elements. 

	The laws defining the variety are as follows: 

	\medskip

	I.  $H(c) \approx  c$,

	$T(S^nT^m(H_k)) \approx  S^nT^{m+1}(H_k)$,

	$S(S^nT^m(H_k)) \approx  S^{n+1}T^m(H_k)$,  

	$S^{-1}(S^nT^m(H_k)) \approx  S^{n-1}T^m(H_k)$,

	$H(S^nT^m(H_k)) \approx  H_{m+k}$ for all $k, m \geq  0$ and $n \in
	\open Z$. 

	\medskip

	II.  $N_Q(q,a,H_k) \approx  \sigma (q,a)$  for all $q \in  Q$, $a 
	\in  \Sigma $, $k \geq  0$,

	$N_H(q,a,H_k) \approx S^{\mu (q,a)}T(H_k)$ for all $q \in Q$, $a \in \Sigma
	$, $k \geq  0$,

	$N_\Sigma (q,a,H_k) \approx  \alpha (q,a)$  for all $q \in  
	Q-Q_{LR}$, $a \in  \Sigma $,  $k \geq  0.$

	\smallskip   

	$N_Q(q_L$, $C_\Sigma (H_k), H_k) \approx  q \approx  N_Q(q_R, 
	C_\Sigma (H_k), H_k)$,

	$N_H(q_L, C_\Sigma (H_k), H_k) \approx  ST(H_k)$,

	$N_H(q_R, C_\Sigma (H_k), H_k) \approx  S^{-1}T(H_k)$,   

	$N_\Sigma (q_L, C_\Sigma (H_k), H_k) \approx  e_L$,

	$N_\Sigma (q_R, C_\Sigma (H_k), H_k) \approx  e_R$, for all $q \in Q
	\setminus Q_{\rm LR}$ and $k \geq 0$. 

	\medskip

	III.  $C_\Sigma T(H_k) \approx N_\Sigma (C_Q(H_k), C_\Sigma (H_k),
	H_k)$,

	$C_QT(H_k) \approx  N_Q(C_Q(H_k)$, $C_\Sigma (H_k), H_k)$,

	$H_{k+1} \approx  N_H(C_Q(H_k), C_\Sigma (H_k), H_k)$ for all $k \geq
	0$. 

	\medskip

	IV.  $C_Q(\lambda ) \approx  C_QH(\lambda )$ for all $\lambda  \in  \Lambda $.

	\medskip

	V.  $R(\lambda ,\lambda ) \approx  0$ for all $\lambda  \in  \Lambda $,   

	$R(\lambda ,S^k(\lambda)) \approx  1$ for all $k > 0$, and $\lambda  \in  
	\Lambda $,

	$F(\lambda ,\lambda ) \approx  0_F$ for all $\lambda  \in  \Lambda $,

	$F(\lambda ,T^k(\lambda) ) \approx 1_F$ for all $\lambda \in \Lambda $
	and $k > 0$, 

	$F(\lambda ,\gamma ) \approx  F(H\lambda , H(\gamma) )$ for all 
	$\lambda ,\gamma  \in  \Lambda $.

	\medskip

	VI.  $C_\Sigma (S^nTH(\lambda) ) \approx  C_\Sigma (S^nH(\lambda) )$ for all 
	$\lambda  \in  \Lambda$ and $n \neq 0$.

	\medskip 

	VII.  $K(0,\lambda ) \approx  \lambda $ for all $\lambda  \in  \Lambda $,

	$K(1,\lambda ) \approx  c$   for all $\lambda  \in  \Lambda $,
	   
	$K(d,d) \approx  0$   for all constants $d$,

	$K(d,e) \approx  1$ for  all constants $d \neq  e \neq  c$,

	$K(\lambda ,d) \approx  1$, for all $\lambda  \in  \Lambda $ and constants $d 
	\neq  c$,

	$K(C_\Sigma( \lambda) ,d) \approx 1$ for all constants $d \notin \Sigma $,
	and $\lambda \in \Lambda $,

	$K(C_Q(\lambda) ,d) \approx  1$  for  all constants $d \notin  Q$,
	and $\lambda  \in  \Lambda $,

	$K(R(\lambda ,\gamma ),d) \approx  1$  for  all $\lambda $, $\gamma  \in  
	\Lambda $, and $d \neq  0,1$,

	$K(F(\lambda ,\gamma ), d) \approx  1$ for all $\lambda $, $\gamma  \in  
	\Lambda $, and $d \neq  0_F, 1_F$.   

	\medskip

	$K(t,t) = 0$ and $K(s,t) = 1$ for all $s,t$ where $s,t$ belong to different 
	members of the following list of sets:  
	$\Lambda $, $\{C_\Sigma (\lambda )|\lambda  \in  \Lambda \}$, 
	$\{C_Q(\lambda )|\lambda  \in  \Lambda \}$, 
	$\{R(\lambda ,\gamma )|\lambda ,\gamma  \in  \Lambda \}$, 
	$\{F(\lambda ,\gamma )|\lambda ,\gamma  \in  \Lambda \}.$ 

	\medskip

	VIII. $EC_Q(\lambda ) \approx  1$   for all $\lambda  \in  \Lambda $   

	$E(h) \approx 0.$ 

	\medskip

	Note that every term generated from $c$ by the operations $S$,
	$S^{-1}$, $T$ and $H$ is equivalent modulo the above laws I to an
	element of $\Lambda $; in fact there is an effective procedure which,
	given such a term $t$, produces $\lambda \in \Lambda $ with $t$
	equivalent (modulo I) to $\lambda $.  Thus we may, and will, ignore
	all such terms $t$ except those in $\Lambda .$

	\medskip

	\noindent
	6.2 NON-UNIFORM SOLVABILITY OF THE WORD PROBLEM

	\medskip

	\noindent
	Proposition 6.1 $V$ {\it does not have uniformly solvable word problem.}

	\medskip

	\noindent
	Proof. The proof is analogous to the proof of Theorem 3.1: for an
	initial tape configuration as described there, we have associated the
	same presentation ${\cal P}$ and prove that the universal Turing
	machine, started on that configuration, eventually halts if and only
	if $0 \equiv _{\cal P} 1.$ 

	As in the proof of Theorem 3.1, the ``only if" part is clear. 

	For the ``if" direction, if the machine does not halt, we again
	produce a model $A \in V$ satisfying all the ${\cal P}$ equations, in
	which $0 \neq 1.$ 

	The underlying set of A is as in the proof of Theorem 3.1, and the operations 
	are as defined there, with the following changes: 

	(ii), (iii) and the definition of $C^*$ are deleted (since we have
	deleted the operations $P$, $U$ and $C^{*}$). 

	In (v)  $R(x,y) = *$ unless both $x,y \in  \Lambda .$ 

	In (vi)  $F(x,y) = *$ unless both $x,y \in  \Lambda .$ 

	In (vii) $K(x,y) = *$ for $x,y$ not in the form of the first two
	lines.  

	In (viii) replace $K^*$ by $K$. 

	\medskip

	\noindent
	6.3 SOLVABILITY OF THE WORD PROBLEM 

	\medskip

	The proof that $V$ has solvable word problem is somewhat different than the 
	proof in section 5.   

	We again differentiate two cases: whether or not ${\cal P}$ has
	degenerate space-time, which in this case means $\lambda \equiv
	_{\cal P} c$ for all $\lambda \in \Lambda .$ 

	In the degenerate case, we have $(\lambda ,c) \in \equiv _{\cal P}$
	for all $\lambda \in \Lambda $, and hence the equations defining our
	variety are equivalent (modulo $\equiv _{\cal P})$ to finitely many
	equations, namely $S(c) = T(c) = S^{-1}(c) = c$ together with all
	instances of the equations defining the variety with $c$ substituted
	for the arbitrary $\lambda \in \Lambda $ which appear.  Thus $\equiv
	_{\cal P}$ is finitely generated relative to the variety of all
	algebras, and hence is decidable by Corollary~1.2. 

	\medskip

	In the non-degenerate case we proceed, at first, similarly to section
	5, bearing in mind that there are fewer space-time elements (see
	above definition of $\Lambda $), but time coordinates, time prefixes,
	and space prefixes are defined as before, except that all these
	notions are always relative to $c = H(c)$, i.e. the only space-time
	component is that of $c$. 

	\medskip

	The proof of the next lemma is essentially the same as the proofs of
	Lemmas 5.2 and 5.3.

	\medskip

	\noindent
	Lemma 6.2 {\it For any finite subset $F \subseteq \Lambda $, with
	maximum time coordinate $m$}, {\it there is a finite $\bar F \subseteq
	\Lambda $ with the same maximum time coordinate such that} 

	(i)  $F \subseteq \overline{F}$  

	(ii) {\it if $\lambda $ is a right subterm of $\gamma \in \bar F$ then
	$\lambda \in  \bar F$}  

	(iii) {\it if  $\lambda  \in  \Lambda $  and  $\lambda  \theta _{\cal P} 
	\gamma $ for $\gamma  \in  \bar F$ then $\lambda  \in  \bar F$}  

	(iv)  {\it if $\lambda  \in  \Lambda $ and $T\lambda  \in  \bar F$ then 
	$\lambda  \in  \bar F.$} 

	\medskip

	Definition of $A$:  

	$A$ consists of all terms appearing in the equations defining the
	variety, and all subterms thereof, plus all elements $C_\Sigma
	(\lambda )$, and $K(\lambda ,\gamma )$ for $\lambda $, $\gamma \in
	\Lambda $, (the elements $C_Q(\lambda )$ are already included).  

	Let  $\equiv _A$ be the restriction to  $A$  of the congruence defining our 
	variety, then, because we have no non-trivial information about  $C_\Sigma $ 
	and $C_Q$, rules III cannot be applied in a non-trivial way, and hence 
	$\equiv _A$ is decidable, as in condition (3) of Proposition 1.1.  Thus $(A, 
	\equiv _A)$ is a partial subalgebra satisfying the hypotheses of Proposition 
	1.1.

	Definition of $B$:  

	Let  $B_{\cal P}$ consist of all terms appearing in the presentation 
	${\cal P}$, and all subterms thereof, and all constants.  

	Enlarge  $B_{\cal P}$ as follows:  

	(i)  For each  $b \in  B_{\cal P}$, if there exists  $a \in  A$ with $a 
	\theta _{\cal P} b$, add one such $a$, and choose $a \in  \Lambda $ whenever 
	possible.  

	(ii)  Let  $F$  consist of all the elements of $\Lambda $ we have so far, and 
	add the set  $\bar F$ of Lemma 6.2.  

	(iii) For all $\lambda ,\gamma  \in  \bar F$, add  $C_Q(\lambda )$, 
	$C_\Sigma (\lambda )$, $R(\lambda ,\gamma )$, $F(\lambda ,\gamma ).$  

	\medskip

	The resulting set is $B$.  It is closed under taking subterms, and for 
	$\lambda ,\gamma  \in  \Lambda $, if $\gamma  \in  B$ and $\lambda  
	\theta _{\cal P} \gamma $ then $\lambda  \in  B.$  

	Let  $\equiv _B$ be  $\theta _{\cal P}|B$, then $B$ and $\equiv _B$ are both 
	finite, and hence decidable. 

	\medskip

	Definition of $A_0$:  

	Let  $A_0 = A \cup  B$, and let  $\equiv _0$ be the partial congruence on $A 
	\cup  B$ generated by $\equiv _A \cup  \equiv _B$; we are going to show that 
	$(A_0, \equiv _0)$ is a partial subalgebra satisfying the hypothesis of 
	Proposition 1.  Since $\theta _{\cal P}$ is generated by  $\equiv _A \cup  
	\equiv _B$ and hence by $\equiv _0$, the decidability of  $\theta _{\cal P}$ 
	will then follow from Proposition 1.1, the proof of which is deferred to the 
	next subsection.  

	Since membership in $A$ is decidable, and $B$ is finite, we know that
	membership in $A_0$ is decidable.  

	Next, we need to establish that  $\equiv _0$ is decidable; this, however, can 
	be proved analogously to the proof in section 5 that $\equiv _0$ (as defined 
	there) is decidable, deleting from that proof consideration of elements which 
	we do not have in this example, such as elements in the image of $U$ or $C^*$ 
	and space-time elements except those in the presently defined $\Lambda $.  

	It remains to check that $(A_0, \equiv _0)$ satisfies hypothesis (3)
	of Proposition 1.1, i.e., that there is an algorithm which, given an
	operation $\sigma $ of arity $n$, and $a_1,\ldots,a_n \in A_0$,
	determines whether there exist $b_1,\ldots,b_n \in A_0$ with $a_i \equiv
	_0 b_i$ and $\sigma (b_1,\ldots,b_n)
	\in  A.$  

	Now, $B$ is finite, and hence we can check all elements of the form 
	$\sigma (b_1,\ldots,b_n) \in  B$, and decide whether  $a_i \equiv _0 b_i$.  

	Thus it is enough to decide whether there exist $b_1,\dots,b_n \in A$
	with $\sigma (b_1,\ldots,b_n) \in A$ and $a_i \equiv _0 b_i$.
	Moreover, if some $a_i \in B - A$ then $a_i \in B_{\cal P}$ and so if
	there exists $b_i \in A$ with $a_i \equiv _0 b_i$ then there exists
	$c_i \in A \cap B$ with $a_i \equiv _0 c_i$ and we may replace $a_i$
	by $c_i$.  Thus we may assume without loss of generality that all
	$a_i \in A.$  

	Hence we have reduced the problem to the following: 
	\begin{quote}
	given
	$a_1,\ldots,a_n \in A$ do there exist $b_1,\ldots,b_n \in A$ with $a_i
	\equiv _0 b_i$ and $\sigma (b_1,\ldots,b_n) \in A$? 
	\end{quote}

	\medskip

	There is an effective procedure which, given $a \in A$, produces $d
	\in A$ such that $d \equiv _A a$ and either $d \in \Lambda $, $d$ is
	a constant, or $d \in C_Q(\Lambda )$, $d \in C_\Sigma (\Lambda )$, $d
	\in R(\Lambda ,\Lambda )$, or $d \in F(\Lambda ,\Lambda )$.  Thus we
	may assume that each $a_i$ is already of this form. 

	We consider the operations in turn  

	$T$: For $a \in A$, if there exists $b \in A$ with $b \equiv _0 a$
	and $T(b) \in A$ then $b \in \Lambda $ hence $a \in \Lambda $.  Thus
	there is such a $b$ if and only if $a \in \Lambda .$ 

	$S,S^{-1}$ and $H$ are the same as $T.$  

	$R$: For $a_1$, $a_2 \in A$, if there exist $b_i \equiv _0 a_i$ with
	$b_i \in A$ and $R(b_1,b_2) \in A$ then $b_i \in \Lambda $ and hence
	$a_i \in 
	\Lambda $, conversely if $a_i \in  \Lambda $ then $R(a_1,a_2) \in  A.$  

	$F$: Is the same as $R$.  

	$C_Q:$ For $b \in A$, $C_Q(b) \in A$ if and only if $b \in \Lambda $,
	hence there exists $b \equiv _0 a$ with $C_Q(b) \in A$ if and only if
	$a \in \Lambda $  

	$C_\Sigma $: same as $C_Q.$  

	$N_Q:$ If $a_1$, $a_2$, $a_3 \in A$ and there exist $b_i \equiv _0
	a_i$ with $N_Q(b_1, b_2, b_3) \in A$, then $b_3 = H_k$ for some $k$,
	then because $a_3 \equiv _0 b_3$ we must have that the time component
	of $a_3$ is $k$ so we can just check whether $a_3 \equiv H_k$.  If
	the answer is affirmative, then for $b_1$ and $b_2$ we have
	$N_Q(b_1, b_2,H_k) \in A$ if and only if either $a_1 \equiv _0
	q \in Q$ and $a_2 \equiv _0 a \in \Sigma $, or $a_1 \equiv _0 q_L$ or
	$q_L$ or $C_Q(H_k)$ and $a_2 \equiv _0 C_\Sigma (H_k)$; there are
	only finitely many cases to check.  

	$N_\Sigma ,N_H$: The argument is the same as for $N_Q.$  

	$E$: For $a \in A$, if there exists $b \in A$ with $E(b) \in A$ then
	either $a \equiv _0 h$ or $a \equiv _0 C_\Sigma (\lambda )$ for some
	$\lambda \in \Lambda $.  The latter occurs if and only if either $a
	\in im(C_\Sigma )$ or $a \equiv _0$ to some $d \in \Sigma $ or $a \equiv
	_0 C_\Sigma (b) \in B$; these finitely many cases can be checked.  

	$K$: If $a_1 \equiv _0 b_1$ and $a_2 \equiv _0 b_2$ and $K(b_1$,
	$b_2) \in A$ then there are various possibilities: 

	(1) $b_2 \in \Lambda $, hence $a_2 \in \Lambda $, and $a_1 \equiv _0
	0$ or $a_1 \equiv _0 1$ or $a_1$ is in $\Lambda $ or in the image of
	$C_\Sigma $, $C_Q$, $R$ or $F$, or $a_1 \equiv _0$ an element in the
	image of one of these operations in $B$.  These finitely many cases
	can be checked. 

	(2) $a_1$ and $a_2$ are each $\equiv _0$ some constant (possibly
	different ones). 

	(3) $a_1 \in \Lambda $ and $a_2$ is either $\equiv _0$ some constant
	not equal to $c$, or $a_2 \in \Lambda $, or $a_2$ is in the image of
	$C_\Sigma $, $C_Q$, $R$ or $F$, or $a_2$ is $\equiv _0$ to an element in
	the image of one of these operations in $B.$ 

	(4) $a_1$ is in the image of $C_\Sigma $ or is $\equiv _0$ an element
	of $B$ which is in the image of $C_\Sigma $, and $a_2 \equiv _0 d$, a
	constant $\notin \Sigma .$ 

	(5) Similar to (4), with $C_\Sigma $ replaced by $C_Q$, or $R$, or
	$F$ respectively, with the appropriate constraint on the constant $d
	\equiv _0 a_2.$

	\noindent
	This completes the proof.

	\medskip

	Applying Proposition~1.1, we have proved the following theorem.

	\medskip

	\noindent
	Theorem 6.4 {\it There is a variety in a finite language defined by a
	recursive set of laws involving only constants which has solvable but
	not uniformly solvable word problem.}

	\medskip

	\noindent
	6.4 PROOF OF PROPOSITION 1.1 

	\medskip

	We complete this section with the promised proof of Proposition 1.1.

	\medskip

	\noindent
	Proof (of Proposition 1.1).  We first produce a partial subalgebra
	$(B, \equiv _B)$ satisfying (1) to (3), such that $A \subseteq B$,
	and $\equiv _A \subseteq \equiv _B$, which has the feature that if
	$a_i \equiv _B b_i$ for $1 \leq i \leq n$ and if $\sigma (a_1,\ldots,
	a_n) \in B$ then $\sigma (b_1,...$, $b_n) \in B.$ 
	$B$ and $\equiv _B$ are defined by induction on the complexity of
	terms, as follows. 

	Let $B_0 = A$ and $\equiv _0$ be $\equiv_A$.
	For each natural number $k$, let 
	\begin{eqnarray*}
	B_{k+1} = B_k &\cup&  \{\sigma (b_1,\ldots, b_n) | \sigma  \in  \Sigma , b_i 
	\in  B_k \mbox{ and  there exists } a_i \equiv _k b_i\\
	              &    & \mbox{ with }\sigma (a_1,\ldots, a_n) \in  A\},
	\end{eqnarray*}

	$R_{k+1} = \{(\sigma (b_1,\ldots, b_n),a) | \sigma \in \Sigma ,b_i
	\in B_k, a \in A \mbox{ and there exists } a_i \equiv _k b_i \mbox{
	with } a = \sigma (a_1,\ldots,a_n)\}.$ 

	Let $\equiv _{k+1} = \equiv _A \cup (\equiv _A \circ R_{k+1}) \cup
	(R^{-1}_{k+1} \circ \equiv _A) \cup (R^{-1}_{k+1} \circ \equiv _A
	\circ R_{k+1}).$ 

	\medskip

	Define $B = \cup  B_k (k \in  \omega )$  and define
	$\equiv _B = \cup  \equiv _k (k \in  \omega )$.

	Note that each $B_k$ is closed under subterms. 

	We will prove the following by induction on $k$: 

	(i)  $\equiv _k \subseteq  \equiv _{k+1}.$ 

	(ii)  If $b \in  B_k$ and $b \equiv _k a, b \equiv _k c$ for $a, c \in  A$ 
	then $a \equiv _A c.$ 

	(iii) $\equiv _k | B_i = \equiv _i$ for all $i < k.$ 

	(iv)  $\equiv _k$ is transitive. 

	(v)  $\equiv _k$ is a partial congruence on $B_k.$ 

	$k = 0$:  trivial. 

	Induction Step:  Suppose we have (i) to (v) for $k$. 

	(i) Then $R_{k+1} \subseteq  R_{k+2}$ and hence $\equiv _{k+1} \subseteq  
	\equiv _{k+2}.$ 

	(ii) Suppose $b \in B_{k+1}$ and $b \equiv _{k+1} a$, $b \equiv
	_{k+1} c$ for $a, c \in A$.  Note that

	\noindent
	$R_{k+1}|A \subseteq \equiv _A$, and hence if $b \in A$ then $b
	\equiv _A a$ and $b \equiv _A c$ so $a \equiv _A c.$ 

	Assume $b \notin A$.  Then $b = \sigma (b_1,\ldots,b_n)$ and there
	exist $a_i \equiv _k b_i$ with

	\noindent
	$\sigma (a_1,\ldots,a_n) \in A$ and $\sigma (a_1,\ldots,a_n) \equiv
	_A a$.  Similarly there exist $c_i \equiv _k b_i$ with $\sigma
	(c_1,...,c_n) \in A$ and $\sigma (c_1,...,c_n) \equiv _A c$.  But by
	the induction hypothesis we get $a_i \equiv _A c_i$ and hence $\sigma
	(a_1,...,a_n) \equiv _A 
	\sigma (c_1,...,c_n)$ so $a \equiv _A c.$ 

	(iii) It is enough to prove that $\equiv _{k+1} | B_k = \equiv _k$ , and for 
	this it is
	enough to prove that $R_{k+1} | B_k \subseteq R_k (or \equiv _A$ if
	$k = 0)$. However, if $(\sigma (b_1,\ldots,b_n), a) \in R_{k+1}$ and
	$b = \sigma (b_1,\ldots,b_n) \in B_k$ then there exist $a_1,\ldots,a_n \in
	A$ with $b_i \equiv _k a_i$ and $a = \sigma (a_1,\ldots,a_n)$.  If $k =
	0$ then we have $b \in A$ and hence $(b, a) \in \equiv _A$.  If $k
	> 0$ then $b \in B_K$ implies that $b_1,\ldots,b_n \in B_{k-1}$ and so the
	induction hypothesis yields $b_i 
	\equiv _{k-1} a_i$ and hence $(b, a) \in  \equiv _k.$ 

	(iv) is an direct consequence of (ii). 

	(v) For $k = 0$ this is just the hypothesis on $(A, \equiv _A)$,
	since we have assumed that $\equiv_A$ is a partial congruence on $A$.
	Suppose $b_i \equiv _{k+1} d_i$
	and $\sigma (b_1,\ldots,b_n) \in B_{k+1}$ and $\sigma
	(d_1,\ldots,d_n) \in B_{k+1}$.  Then $b_i$, $d_i \in B_k$ for $1 \leq
	i \leq n$ and hence by (iii), $b_i \equiv _k d_i$.  Also, there exist
	$a_i$, $c_i \in A$, $1 \leq i \leq n$ with $b_i \equiv _k a_i$ and
	$d_i \equiv _k c_i$ and $\sigma (a_1,\ldots,a_n)$, $\sigma (c_1,\ldots,c_n)
	\in A$.  By (iv) and (ii) we obtain $a_i \equiv _A c_i$ and hence
	$\sigma (a_1,\dots,a_n) \equiv _A 
	\sigma (c_1,\dots,c_n)$ and thus $\sigma (b_1,\ldots,b_n) \equiv _{k+1} 
	\sigma (d_1,\ldots,d_n)$. 

	It remains to verify (1), (2), and (3) for $(B, \equiv _B).$ 

	Note first of all that $\sigma (b_1,\ldots,b_n) \in B$ if and only if
	there exist $a_1,\ldots,a_n \in A$ with $a_i \equiv _B b_i$ and
	$\sigma (a_1,\ldots,a_n) \in A$.  Moreover, if $\sigma (b_1,\dots,b_n) \in
	B$ and the $a_i$ are as above then whenever $b_i \equiv _B c_i \in A$
	with $\sigma (c_1,\ldots,c_n) \in A$ then $\sigma (a_1,\ldots,a_n) \equiv
	_A \sigma (c_1,\ldots,c_n).$ 

	We show that there is an algorithm which, given $b \in FX$,
	determines whether $b \in B$ and in the affirmative case produces $a
	\in A$ with $b \equiv _B a$. 

	Consider $b \in FX$.  It is decidable whether $b \in A$, and in the
	affirmative case we are finished.  If $b \notin A$ then $b = \sigma
	(b_1,\dots,b_n)$ for unique $b_1,\ldots,b_n$ and $\sigma $.  In this
	case, $b \in B$ if and only if all the $b_i \in B$, and there exist
	$a_i \in A$ with $\sigma (a_1,\dots,a_n) \in A$ and $b_i \equiv _B
	a_i$.  The $b_i$ are of lower complexity than $b$; determine for each
	whether it belongs to $B$ and if so, produce $c_i \in A$ with $c_i
	\equiv _B b_i$.  Given the $c_i$ it is decidable whether there exist
	$a_1,\ldots,a_n \in A$ with $a_i \equiv _A c_i$ and $\sigma
	(a_1,\ldots,a_n) \in A$, and moreover, since membership in $A$ is
	decidable, we can effectively produce the $a_i$ in the affirmative
	case, thus yielding an appropriate $a \in A$ with $a \equiv _B b$,
	namely $a = \sigma (a_1,\ldots,a_n) \equiv_n \sigma(c_1, \ldots, c_n)
	\equiv_n \sigma(b_1, \ldots, b_n) = b.$ 

	Thus membership in $B$ is decidable, and hence $\equiv _B$ is
	decidable: given $b$, $c \in B$ we effectively produce $a, d \in A$
	with $b \equiv _B a$, $c \equiv _B d$ and then $b \equiv _B c$ if and
	only if $a \equiv _A d$, and the latter is decidable. 

	Finally, given an $n$-ary operation $\sigma $ and $b_1,\dots,b_n \in
	B$, there exist $a_1,\ldots,a_n \in B$ with $a_i \equiv _B b_i$ and
	$\sigma (a_1,\ldots,a_n) \in B$, if and only if $\sigma (b_1,\ldots,b_n)
	\in B$, and we have just proved that the latter is decidable. 

	Thus $(B, \equiv _B)$ is a partial subalgebra with all the
	properties claimed above. 

	Now, each element $s \in FX$ can effectively be written as $s =
	s'(u_1,\ldots,u_n)$ where the $u_i \in B$ are subterms which are
	maximal with respect to belonging to $B.$ 

	Define a relation $\equiv$ on $FX$ as follows: for $s =
	s'(u_1,\ldots,u_n)$ and $t = t'(v_1,\ldots,v_k)$, where the $u_i$, $v_j$
	are maximal $B$-subterms, 

	$s \equiv  t$ if and only if $s' = t'$ and $u_i \equiv _B v_i$ for all $i.$

	\noindent
	Then $\equiv$  is a congruence on $FX$, which extends $\equiv _A$ and is 
	generated by it, and so is the congruence on $FX$ generated by $\equiv _A$.  
	Moreover, the above description of $\equiv $, together with the decidability 
	of $\equiv _B$, yields the decidability of $\equiv $, as required.

	\bigskip

	\noindent
	{\S}7. A VARIETY WITH INFINITELY MANY OPERATIONS 

	\medskip

	If we allow infinitely many operations, then it is much easier to
	obtain a variety with solvable, but not uniformly solvable, word
	problem.  The following example is a modification of an example given
	in Wells [W, p.161] for a different, although related purpose. He
	suggested its relevance to our question. 

	\medskip

	Let $V$ be the variety with a constant, 0, a binary operation denoted by 
	juxtaposition, and countably many unary operations  $h_n$ $(n \in  \omega )$ 
	satisfying the following laws:    

	\medskip

	$xy \approx  yx$    

	$x(yz) \approx  (xy)z$    

	$x0 \approx  0$    

	$x^2 \approx  0$    

	$xh_n(y) \approx  0$           for all $n \in  \omega $    

	$h_n(h_n(x)) \approx  h_n(x)$    for all $n \in  \omega $    

	$h_n(h_k(x)) \approx  0$        for all $n \neq  k$

	\medskip

	\noindent
	and 

	$(*)$    ${h_{m_n}}^n(x_1x_2  \cdots x_{m_n}) \approx  0.$ 

	\medskip

	where  $\{m_n | n \in  {\open N}\}$ is a recursive listing of a non-recursive 
	set $X$. 

	\medskip

	Thus $V$ is a variety of commutative, square-zero semigroups with
	countably many idempotent unary operations, and the above is a
	recursive set of equations defining $V$. It is worth commenting on why
	the system of equations is recursive.  Obviously the only problem is
	identifying when an equation is included in the scheme $(*)$. Now the
	equations in $(*)$ are of the form ${h_j}^k(x_1\cdots x_j) \approx 0$.
	Such an equation is in $(*)$ if and only if $j = m_k$. The trick of
	using ${h_{m_n}}^n(x_1 \cdots x_{m_n})$ rather than $h_{m_n}(x_1\cdots
	x_{m_n})$ is a variant of the old trick of pleonasm due to Craig
	 which he used to prove that any theory with a
	recursively enumerable axiomatization has a recursive axiomatization
	(see Monk [7, p.262]) . 

	\medskip

	We will show that $V$ has an undecidable equational theory, and hence
	does not have uniformly solvable word problem, by establishing that
	$V$ satisfies the equation $h_k(x_1 \cdots x_k) \approx 0$ if and only if
	$k \in X$. One direction is trivial by the laws above. To complete
	the proof of undecidibility, we construct an algebra in which for $k
	\notin X$, $h_k$ is non-zero on a product of $k$ elements.

	Let $S$ be the free algebra on countably many generators in the class
	of commutative semigroups with 0 satisfying $x^2 = 0$.  Let $\{a_i |
	i \in {\open N}\}$ be a countable set disjoint from $S$, and let $A =
	S \cup \{a_i | i \in {\open N}\}$. Define the operations in $A$ as
	follows: the binary multiplication extends that of $S$, and otherwise
	is constant with value 0. 

	For $n \in  X$, $h_n$ is constant with value 0. For $n \notin  X$, 
	$$h_n(x) 
	=\left\{ \begin{array}{ll}
	0 &\hbox{  if } x = 0\hbox{  or } x = a_i\hbox{ for some } i  \neq  n \\
	a_n & \hbox{ otherwise}
	\end{array}
	\right. .$$ 

	\noindent
	It is easy to check that this algebra has the desired properties. 

	\medskip

	Now, to see that $V$ has solvable word problem, consider a finite
	presentation ${\cal P}$ in generators $b_1,\ldots,b_n$.  Let $m$ be
	greater than
	$n$, and greater than $k$ for any $k$ such that $h_k$ appears in one of the
	defining relations of ${\cal P}$.  Let $B$ be the algebra given by
	the presentation ${\cal P}$ in the variety $V'$, which has operations
	0, multiplication, and $h_i$ for $i \leq m$, and is defined by the
	laws defining $V$ which involve only the $h_i$ for $i \leq m$. Then
	$B$ is finite, and hence the word problem for ${\cal P}$ relative to
	the variety $V'$ is decidable.  Let $C \subseteq B$ consist of all
	non-zero elements of $B$ which are not the image of any $h_i$ $(i \leq
	m)$.  Then the algebra $A$ given by the presentation ${\cal P}$ in the
	variety $V$ has as underlying set $B \cup (C \times \{i \in {\open
	N}|i > m\})$; the multiplication extends that of $B$ and otherwise
	has value 0, the $h_i$ for $i \leq m$ extend those of $B$ and
	otherwise have value 0, and for $i > m$ and $c \in C$, $h_i(c) =
	(c,i) = h_i((c,i))$ and $h_i$ has value 0 otherwise. The equations
	for $i > m$ are satisfied in $A$ because all products $x_1\cdots x_k$ for
	$k > m$ are 0. (It is a simple exercise to show that the laws imply
	$h_k(0) \approx 0$.) This explicit description of $A$ yields a solution to
	the word problem for ${\cal P}$ relative to the variety $V.$

	\medskip

	\noindent
	REFERENCES.

	\noindent
	[1] S. Burris and H. Sankapannavar.  A Course in Universal Algebra 
	Springer-Verlag, 1981.

	\noindent
	[2] T. Evans. The Word Problem for Abstract Algebras.  J. London
	Math. Society,   26 (1951) 64--71.

	\noindent
	[3] T.Evans. Embeddability and the Word Problem.  J. London Math.
	Society  28 (1953) 76--80.

	\noindent
	[4] T. Evans.  Some Solvable Word Problems. in Word Problems II ed.
	by  S. I. Adian, W. W. Boone and G. Higman.  North Holland (1980) 87--100.

	\noindent
	[5] G. Hutchinson.  Recursively unsolvable word problems of modular
	lattices  and diagram chasing. J. Alg. 26 (1973) 385--399.

	\noindent
	[6] A. Markov.  On the impossibility of certain algorithms in the
	theory  of associative systems. Dokl. Akad. Nauk. SSSR (NS) 55 (1947)
	583--586. 

	\noindent
	[7] D. Monk.  Mathematical Logic. Springer-Verlag 1976. 

	\noindent
	[8] P.S. Novikov.  On the Algorithmic unsolvability of the word problem in  
	group theory.  Trudy Mat. Inst. Steklov 44, English Translation Proc. 
	Steklov Inst. Math. (2) 9 (1958) 1--122.

	\noindent
	[9] E. Post.  Recursive unsolvability of a problem of Thue.  Journal
	of Symbolic Logic. 12 (1947) 1--11.

	\noindent
	[10] B. L. van der Waerden.  Geometry and Algebra in Ancient
	Civilizations.  Springer-Verlag 1983.

	\noindent
	[11] B. Wells.  Pseudorecursive varieties and their implications for
	word problems Ph.D. thesis, University of California, Berkeley 1982.

\end{document}